\documentclass[a4paper,10pt]{amsproc}
\title{Rings of functions whose closure of discontinuity set is in an ideal of closed sets}
\usepackage{amsfonts, amsmath, amssymb, graphicx, mathrsfs}
\usepackage[all]{xy}

\newdir{ >}{!/6pt/@{ }*:(1,0.2)@_{>}*:(1,-0.2)@^{>}}

\theoremstyle{plain}

\newtheorem{theorem}{Theorem}[section]
\newtheorem{lemma}[theorem]{Lemma}
\newtheorem{proposition}[theorem]{Proposition}
\theoremstyle{definition}
\newtheorem{definition}[theorem]{Definition}
\theoremstyle{definitions}
\newtheorem{definitions}[theorem]{Definitions}
\theoremstyle{observation}
\newtheorem{observation}[theorem]{Observation}

\newtheorem{counter example}[theorem]{Counter Example}
\newtheorem{notation}[theorem]{Notation}

\newtheorem{corollary}[theorem]{Corollary}
\newtheorem{example}[theorem]{Example}

\numberwithin{equation}{section}

\begin{document}

\author[A. Dey]{Amrita Dey}	\address{Department of Pure Mathematics, University of Calcutta, 35, Ballygunge Circular Road, Kolkata 700019, West Bengal, India}	\email{deyamrita0123@gmail.com}

\author[S. K. Acharyya]{Sudip Kumar Acharyya} \address{Department of Pure Mathematics, University of Calcutta, 35, Ballygunge Circular Road, Kolkata 700019, West Bengal, India}  \email{sdpacharyya@gmail.com}

\author[S. Bag]{Sagarmoy Bag}	\address{Department of Mathematics, Bangabasi Evening College, 19, Raj Kumar Chakraborty Ln.,Baithakkhana, Kolkata  West Bengal, India}	\email{sagarmoy.bag01@gmail.com}

\author[D. Mandal]{Dhananjoy Mandal} \address{Department of Pure Mathematics, University of Calcutta, 35, Ballygunge Circular Road, Kolkata 700019, West Bengal, India}  \email{dmandal.cu@gmail.com}



	\title[On $C(X)_\mathcal{P}$]{Rings of functions whose closure of discontinuity set is in an ideal of closed sets}

	%

	\thanks {}
	
	\maketitle
	
	\large
	\begin{abstract}
		
		Let $\mathcal{P}$ be an ideal of closed subsets of a topological space $X$. Consider the ring, $C(X)_\mathcal{P}$ of real valued functions on $X$ whose closure of discontinuity set is a member of $\mathcal{P}$. We investigate the ring properties of $C(X)_\mathcal{P}$ for different choices of $\mathcal{P}$, such as the $\aleph_0$-self injectivity and regularity of the ring, if and when the ring is Artinian and/or Noetherian. The concept of $\mathcal{F}P$-space was introduced by Z. Gharabaghi, M. Ghirati and A. Taherifar in 2018 in a paper published in Houston Journal of Mathematics. In this paper, they established a result stating that every $P$-space is a $\mathcal{F}P$-space. We furnish that this theorem might fail if $X$ is not Tychonoff and we provide a suitable counter example to prove our assertion.
		
	\end{abstract}

	\section{Introduction}
	
	We assume all our topological spaces $(X,\tau)$ to be $T_1$.  Suppose, $\mathcal{P}$ is an ideal of closed subsets of the topological space $(X,\tau)$, that is, $\mathcal{P}$ is a collection of closed subsets of $X$ satisfying the following conditions : 
	\begin{enumerate}
		\item $A,B\in \mathcal{P}\implies A\cup B \in \mathcal{P}$, and
		\item If $A\in \mathcal{P}$ and $B\subseteq A$ with $B$ closed in $X$, then $B\in \mathcal{P}$.
	\end{enumerate} For example, the set of all closed subsets of a topological space $(X,\tau)$ is an ideal of closed subsets of $X$. Other examples include the set of all finite subsets of $X$, denoted by $\mathcal{P}_f$; the set of all closed compact subsets of $X$, denoted by $\mathcal{K}$ and even the set of all closed Lindelof subsets of $X$, denoted by $\mathcal{L}$. Another trivial but important example of an ideal of closed sets is $\{\emptyset \}$.
	We introduce the triplet $(X,\tau,\mathcal{P})$ and call it a $\tau \mathcal{P}$-space. For a subset $S$ of $X$, $\mathcal{P}_S=\{P\cap S\colon P\in \mathcal{P} \}$ is an ideal of closed subsets of $S$ and we say $(S,\tau_S,\mathcal{P}_S)$ is a $\tau \mathcal{P}$-subspace of $X$, where $\tau_S$ is the subspace topology on $S$ induced by the topology $\tau$ on $X$. 
	
	We now define the set $C(X)_\mathcal{P}=\{f\in \mathbb{R}^X\colon \overline{D_f}\in \mathcal{P} \}$, where $D_f$ denotes the set of all points of discontinuity of $f$. It can be easily seen that $C(X)_\mathcal{P}$ forms a ring with respect to pointwise addition and multiplication. Further, on defining $(f\wedge g)(x)=\min\{f(x),g(x)\}, \forall x\in X$ and $(f\vee g)(x)=\max\{f(x),g(x)\}, \forall x\in X$, for any $f,g\in C(X)_\mathcal{P}$, $(C(X)_\mathcal{P},+,\cdot,\vee,\wedge)$ forms a lattice ordered ring. Also, the set $C^*(X)_\mathcal{P}$ of all bounded functions in $C(X)_\mathcal{P}$ is a lattice ordered subring of $C(X)_\mathcal{P}$. Throughout this article, we say $X$ is a $\tau\mathcal{P}$-space, instead of $(X,\tau,\mathcal{P})$ and $C(X)_\mathcal{P}$ instead of $(C(X)_\mathcal{P},+,\cdot,\vee,\wedge)$. It is clear that for $\mathcal{P}=\{\emptyset \}$, $C(X)_\mathcal{P}=C(X)$, the ring of real valued continuous functions, which has been studied thoroughly in \cite{GJ1976} and when $\mathcal{P}$ is the set of all closed subsets of $X$, $C(X)_\mathcal{P}=\mathbb{R}^X$. Further, when $\mathcal{P}$ is the set of all finite subsets of $X$, then $C(X)_\mathcal{P}=C(X)_F$, the ring of functions on $X$ having finitely many discontinuities \cite{GGT2018} and when $\mathcal{P}$ is the set of all closed nowhere dense subsets of $X$, $C(X)_\mathcal{P}=T'(X)$ (The ring of those real valued continuous functions on $X$ for which there exists a dense open subset $D$ of $X$ such that $f|_D\in C(D)$ \cite{A2010}). We denote the ring $C(X)_\mathcal{P}$ as $C(X)_K$, when $\mathcal{P=K}$. 
	
	The intent of this article is to study and discuss various properties of the ring $C(X)_\mathcal{P}$.
	
	All rings are assumed to be commutative rings with identity. An ideal of a ring is considered to be a proper subset of the ring throughout this article. An ideal of a ring $R$ is said to be essential if it intersects every non-zero ideal of $R$ non-trivially. The sum of all minimal ideals of $R$ is called its socle. $rad(R)$ denotes the Jacobson radical of $R$, that is $rad (R)$ is the intersection of all maximal ideals of $R$. A ring $S$ containing a reduced ring $R$ is called a ring of quotients of $R$ if and only if for each $s\in S\setminus \{0\}$, there exists $r\in R$ such that $sr\in R\setminus \{0\}$. A ring $R$ is said to be almost regular if for every non unit element $a\in R$, there exists a non identity element $b\in R$ such that $ab=a$. We use the notation $\chi_{_A}$ for the characteristic function of $A$ on $X$, defined by \[ \chi_{_A}(x)=\begin{cases}
		1 \text{ when }x\in A, \\ 0 \text{ otherwise}
	\end{cases}.\]
	For an ideal $I$ of $R$, $Ann(I)=\{r\in R\colon ri=0 \text{ for all }i\in I \}$. A ring $R$ is said to be an $IN$-ring if for any two ideals $I$ and $J$, $Ann(I \cap J)=Ann (I) \cap Ann(J)$. A ring $R$ is said to be an $SA$-ring if for any two ideals $I$ and $J$, there exists an ideal $K$ of $R$ such that $Ann (I)\cap Ann(J)=Ann(K)$. Clearly, an $SA$-ring is always an $IN$-ring. A ring $R$ is said to be a Baer ring if annihilator of every non-empty ideal is generated by an idempotent. The next lemma, proved in \cite{T2014} states a characterisation for $IN$-rings for reduced rings.
	
	\begin{lemma} \label{1.1}
		Let $R$ be a reduced ring. Then the following statements are equivalent.
		\begin{enumerate}
			\item For any two orthogonal ideals $I$ and $J$ of $R$, $Ann(I)+Ann(J)=R$.
			\item For any two ideals $I$ and $J$ of $R$, $Ann(I)+Ann(J)=Ann(I\cap J)$.
		\end{enumerate}
	\end{lemma} 
	Further, Birkenmeier, Ghirati and Taherifar established the following set of equivalent conditions in \cite{BGT2015}.
	\begin{lemma} \label{1.2}
		Let $R$ be a reduced commutative ring. Then the following statements are equivalent.
		\begin{enumerate}
			\item $R$ is a Baer ring.
			\item R is a $SA$-ring.
			\item The space of prime ideals of $R$ is extremally disconnected.
			\item $R$ is an $IN$-ring.
		\end{enumerate}
	\end{lemma}
		
	In Section 2, we define zero sets of the form $Z_\mathcal{P}(f)=\{x\in X\colon f(x)=0 \}$, for $f\in C(X)_\mathcal{P}$. We define a cozero set of a function in $C(X)_\mathcal{P}$ to be the complement of $Z_\mathcal{P}(f)$ and denote it by $coz(f)$. We denote the collection of all zero sets of functions in $C(X)_\mathcal{P}$ by $Z_\mathcal{P}[X]$ and the set of all cozero sets of functions in $C(X)_\mathcal{P}$ by $coz[X]$. Also, for a subset $S\subseteq C(X)_\mathcal{P}$, we write $Z_\mathcal{P}[S]=\{Z_\mathcal{P}(f)\colon f\in S \}$ and $coz[S]=\{coz(f)\colon f\in S \}$. Here, we also introduce the notion of $\mathcal{P}$-completely separated subsets of $X$ and achieve its characterisation via zero sets of functions in $C(X)_\mathcal{P}$. Next, we define $z_\mathcal{P}$-filters on $X$ and $z_\mathcal{P}$-ideals in $C(X)_\mathcal{P}$, and examine the duality existing between them. As expected it is  realised that there exists a one-to-one correspondence between the set of all maximal ideals $Max(C(X)_\mathcal{P})$ of $C(X)_\mathcal{P}$ and the set $\mathcal{U}(X)_\mathcal{P}$ of all maximal $z_\mathcal{P}$-filters on $X$. We exploit this duality to show that $M(C(X)_\mathcal{P})$ equipped with the hull-kernel topology is homeomorphic to $\mathcal{U}_\mathcal{P}$ equipped with the Stone topology [Theorem \ref{2.9}].
	
	In Section 3, we examine when does $C(X)_\mathcal{P}$ become closed under uniform limit. When $C(X)_\mathcal{P}$ is closed uniform limit, we then say $X$ is a $\tau\mathcal{PU}$-space. It is clear that when $\mathcal{P}=\{\emptyset \}$, $C(X)_\mathcal{P}=C(X)$, which is closed under uniform limit for any topological space $X$. We also check that for any choice of $\mathcal{P}$, if the set of all non-isolated points of $X$ is a member of $\mathcal{P}$, then  $(X,\tau,\mathcal{P})$ is a $\tau\mathcal{PU}$-space. For some special choice of $\mathcal{P}$, the converse of this result is seen to be valid [Theorem \ref{3.2} and \ref{3.5}]. It is further noted that when $C(X)_\mathcal{P}$ is isomorphic to $C(Y)$, for some topological space $Y$, then $X$ is a $\tau\mathcal{PU}$-space. Using the above results, we have given an alternative proof of Theorem 3.4 in \cite{GGT2018}. At the end of this section, we establish a result analogous to the Urysohn's extension Theorem for $C(X)$, stated in \cite{GJ1976}, for a $\tau\mathcal{PU}$-space.
	
	In Section 4, we define a $\tau\mathcal{P}$-space $X$ to be $\tau\mathcal{P}$-compact if every family of zero sets in $X$, having the finite intersection property has a non-empty intersection. We obtain a characterisation of $\tau\mathcal{P}$-compact spaces using fixed $z_\mathcal{P}$-filters of $X$ and fixed ideals of $C(X)_\mathcal{P}$. Incidentally, we also define $\tau\mathcal{P}$-pseudocompact spaces as follows: a $\tau\mathcal{P}$-space $X$ is $\tau\mathcal{P}$-pseudocompact if and only if $C(X)_\mathcal{P}=C^*(X)_\mathcal{P}$. We go on to show that every $\tau\mathcal{P}$-compact space is $\tau\mathcal{P}$-pseudocompact. In the same section, we call a maximal ideal $M$ of $C(X)_\mathcal{P}$ to be real if $C(X)_\mathcal{P}/M$ is isomorphic to $\mathbb{R}$. We define a $\tau\mathcal{P}$-space to be $\tau\mathcal{P}$-real compact if every real maximal ideal of $C(X)_\mathcal{P}$ is fixed. We establish that a $\tau\mathcal{P}$-space is $\tau\mathcal{P}$-compact if and only if it is both $\tau\mathcal{P}$-pseudocompact and $\tau\mathcal{P}$-real compact. For $\mathcal{P}=\{\emptyset \}$, this reads as follows : a topological space $X$ is compact if and only if it is both pseudocompact and real compact, as proved in \cite{GJ1976}. We construct examples to ensure that a $\tau \mathcal{P}$-pseudocompact space may not be $\tau \mathcal{P}$-compact and a real $\tau \mathcal{P}$-compact space need not be $\tau \mathcal{P}$-compact.
	
	In Section 5, we continue our study of $C(X)_\mathcal{P}$ with the additional hypothesis that each singleton subset of $X$ is a member of $\mathcal{P}$. It follows that $\chi_{\{x\}}\in C(X)_\mathcal{P}$, for all $x\in X$. Under this restriction, we see that $C(X)_\mathcal{P}$ is an almost regular ring and any $f\in C(X)_\mathcal{P}$ is either a zero divisor or an unit. We further show that $C(X)_\mathcal{P}=C(X)$ if and only if $X$ is discrete if and only if $C(X)_\mathcal{P}$ is a ring of quotients of $C(X)$. Further, we are able to establish that a necessary and sufficient condition for a $\tau \mathcal{P}$-space to be a $\tau \mathcal{P}$-compact space is that $X$ is finite. We find out a necessary and sufficient conditions under which an ideal of $C(X)_\mathcal{P}$ is a minimal ideal and establish that the socle of $C(X)_\mathcal{P}$ consists of all functions that vanish everywhere except on a finite set and is itself an essential ideal that is also free. We further note that $Soc(C(X)_\mathcal{P})=C(X)_\mathcal{P}\iff X$ is finite. Exploiting these results, we establish that $C(X)_\mathcal{P}$ is an Artinian (and Noetherian) ring if and only if $X$ is finite. We complete this section by providing a set of conditions equivalent to $C(X)_\mathcal{P}$ being an $IN$-ring, $SA$-ring and Baer ring. We have also provided counter examples to show that these results may fail when the restriction "$\mathcal{P}$ contains all singleton subsets of $X$" is lifted.
	
	In Section 6, we examine the regularity of $C(X)_\mathcal{P}$. Here, we define a $\tau\mathcal{P}$-space, ($X,\tau,\mathcal{P}$) to be a $\mathcal{P}P$-space if $C(X)_\mathcal{P}$ is regular. We show that a $P$-space is a $\mathcal{P}P$-space, when $X$ is Tychonoff. We further provide a counter example to show that this might fail when $X$ is not Tychonoff. This counter example also shows that Theorem 6.1 in \cite{GGT2018} fails when $X$ is not Tychonoff. We conclude this section by giving a characterisation of a $\mathcal{P}P$-space, using the members of $\mathcal{P}$.
	
	Finally, in the seventh section, we use the concept of $\phi$-algebra and an algebra of measurable functions to establish a condition involving a $\tau \mathcal{PU}$-space, under which $C(X)_\mathcal{P}$ is $\aleph_0$-self injective. We also provide an example that shows that the condition $X$ is a $\tau \mathcal{PU}$-space is not superfluous

	\section{Definitions and Preliminaries}
	
	\begin{notation}
		Let $\mathcal{P}'$ be the ideal of all closed subsets of the set of isolated points of $X$.
	\end{notation}

	\begin{theorem} \label{2.03}
		$C(X)_\mathcal{P}=C(X)$ if and only if $\mathcal{P}\subseteq \mathcal{P'}$.
	\end{theorem}
	\begin{proof}
		Let $\mathcal{P}\nsubseteq \mathcal{P'}$. Then there exists $A\in \mathcal{P}$ such that there exists $x_0\in A$ which is a non-isolated point in $X$. Let $f=\chi_{\{x_0 \}}$. Then $f\in C(X)_\mathcal{P}$. However, since $x_0$ is a non-isolated point of $X$, $f\notin C(X)$.
		The converse is obvious.
	\end{proof}
	
	We have stated in the introduction that for $\mathcal{P}=\mathcal{P}_{nd}$, $C(X)_\mathcal{P}=T'(X)$, where $T'(X)$ is the ring of those real valued continuous functions on $X$ for which there exists a dense open subset $D$ of $X$ such that $f|_D\in C(D)$ \cite{A1997}. We give a proof supporting this statement.
	
	\begin{theorem} \label{2.05}
		$C(X)_{\mathcal{P}_{nd}}=T'(X)$ 
	\end{theorem}
	\begin{proof}
		Let $f\in C(X)_{\mathcal{P}_{nd}}$, then $\overline{D_f}$ is nowhere dense. Therefore $X\setminus \overline{D_f}$ is dense in $X$. Also, $\overline{D_f}$ is closed in $X$ $\implies$ $X\setminus \overline{D_f}$ is open in $X$. Finally, $f$ is continuous on $X\setminus D_f \supseteq X\setminus \overline{D_f}$. Thus $f\in T'(X)$.
		Conversely, let $g\in T'(X)$. Then there exists an open dense subset $D$ of $X$ such that $f|_D\in C(D)$ $\implies$ $D_f\subseteq X\setminus D.$ So $\overline{D_f}\subseteq \overline{X\setminus D}=X\setminus D$. Since the complement of an open dense set is a closed nowhere dense set, $X\setminus D\in \mathcal{P}_{nd}$ and so $\overline{D_f}\in \mathcal{P}_{nd}$. Therefore $f\in C(X)_{\mathcal{P}_{nd}}$.
	\end{proof}
	
		We define $\mathcal{I}$ to be the family of all ideals of closed subsets of $X$ and $\mathcal{J}(X)$ to be the family of all subrings of $\mathbb{R}^X$ containing $C(X)$. Both of these families form a lattice with respect to the subset inclusion. 
	
	It can be easily observed that for two rings $S_1,S_2\in \mathcal{J}(X)$, \\$\displaystyle{S_1\vee S_2=\{\sum_{i=1}^{m}f_ig_i\colon f_1\in S_1, g_i\in S_2, i=1,2,...,m. \}}$ is the smallest subring of $\mathbb{R}^X$ containing $S_1\cup S_2$ and for two ideals of closed sets $\mathcal{P} $ and $\mathcal{Q}$ in $\mathcal{I}$, $\mathcal{P}\vee \mathcal{Q}=\{A\cup B\colon A\in \mathcal{P} \text{ and }B\in \mathcal{Q} \}$ is the smallest ideal of closed subsets of $X$ containing $\mathcal{P}$ and $\mathcal{{Q}}$. It is obvious that for two rings $S_1,S_2\in \mathcal{J}(X)$, $\displaystyle{S_1\wedge S_2=S_1\cap S_2}$ and for two ideals of closed sets $\mathcal{P} $ and $\mathcal{Q}$ in $\mathcal{I}$, $\mathcal{P}\wedge \mathcal{Q}=\mathcal{P}\cap \mathcal{Q}$. So $C(X)_{\mathcal{P}\wedge \mathcal{Q}}=C(X)_\mathcal{P}\wedge C(X)_\mathcal{Q}$. However, $C(X)_{\mathcal{P}\vee \mathcal{Q}}=C(X)_\mathcal{P}\vee C(X)_\mathcal{Q}$ may not hold for all $T_1$-spaces but is true for a topological space $X$ if every open subspace of $X$ is $C$-embedded.
	
	\begin{theorem} \label{lattice}
		Let $X$ be a topological space where every open subspace of $X$ is $C$-embedded. Then for any $\mathcal{P}$ and $\mathcal{Q}$ in $\mathcal{I}$, $C(X)_{\mathcal{P}\vee \mathcal{Q}}=C(X)_\mathcal{P}\vee C(X)_\mathcal{Q}$.
	\end{theorem}
	\begin{proof}
		Let $\displaystyle{\alpha=\sum_{i=1}^{m}f_ig_i\in C(X)_\mathcal{P}\vee C(X)_\mathcal{Q}}$. Then\\ $\displaystyle{D_{\alpha}\subseteq \bigcup_{i=1}^{m} D_{f_ig_i}\subseteq \bigcup_{i=1}^{m} (D_{f_i}\cup D_{g_i})}$ where $\overline{D_{f_i}}\in \mathcal{P}$ and $\overline{D_{g_i}}\in \mathcal{Q}$ for all $i=1,2,...,m$. So $\overline{D_{f_i}\cup D_{g_i}}\in \mathcal{P}\vee \mathcal{Q}$ for all $i=1,2,...,m$. Since $\mathcal{P}\vee \mathcal{Q}$ is an ideal of closed sets, $\overline{D_\alpha}\in \mathcal{P}\vee \mathcal{Q}$ and hence $\alpha\in C(X)_{\mathcal{P}\vee \mathcal{Q}}$. Conversely, let $f\in C(X)_{\mathcal{P}\vee \mathcal{Q}}$. Then $f|_{X\setminus \overline{D_f}}$ is continuous on the open subspace $X\setminus \overline{D_f}$ of $X$. By our hypothesis there exists $\widehat{f}\in C(X)$ such that $\widehat{f}|_{X\setminus \overline{D_f}}=f|_{X\setminus \overline{D_f}}$. Also $\overline{D_f}\in \mathcal{P}\vee \mathcal{Q}$. So $\overline{D_f}=A\cup B$ where $A\in \mathcal{P}$ and $B\in \mathcal{Q}$. We define $g\colon X\longrightarrow \mathbb{R}$ by $g(x)=\begin{cases}
			\widehat{f}(x) \text{ when }x\in X\setminus A \\ 
			f(x) \text{ when }x\in A\setminus B \\
			\frac{1}{2}f(x) \text{ when }x\in A\cap B
		\end{cases}$ and $h\colon X\longrightarrow \mathbb{R}$ by $h(x)=\begin{cases}
			0 \text{ when }x\in X\setminus B\\ 
			f(x)-\widehat{f}(x) \text{ when }x\in B\setminus A \\
			\frac{1}{2}f(x) \text{ when }x\in A\cap B
		\end{cases}.$  Then $\overline{D_g}\subseteq A$ and $\overline{D_h}\subseteq B$. Therefore $g\in C(X)_\mathcal{P}$ and $h\in C(X)_\mathcal{Q}$ and $f=g+h\in C(X)_\mathcal{P}\vee C(X)_\mathcal{Q}$.
	\end{proof}
	It is important to note that if $X$ is an irreducible spaces, then every open subspace of $X$ is $C$-embedded and thus for an irreducible space $X$,  $C(X)_{\mathcal{P}\vee \mathcal{Q}}=C(X)_\mathcal{P}\vee C(X)_\mathcal{Q}$.
	
	We summarise all this in the following theorem.
	\begin{theorem} \label{lattice2}
		Let $X$ be a topological space where every open subspace of $X$ is $C$-embedded. Then $\phi \colon \mathcal{I}\longrightarrow \mathcal{J}(X)$ defined by $\phi(\mathcal{P})=C(X)_\mathcal{P}$ is a lattice homomorphism.
	\end{theorem}

	We recall that $f \in C(X)_\mathcal{P}$, $Z_\mathcal{P}(f) = \{x\in X \colon f(x)=0\}$ is called a zero set of $f$.
	Also, a subset $A$ of $X$ is called a zero set if $A=Z_\mathcal{P}(f)$, for some $f\in C(X)_\mathcal{P}$. 
	Let $Z_\mathcal{P}[X]$ be the set of all zero sets of $X$. 
	
	It is easy to verify that:
	\begin{enumerate}
		\item For $f\in C(X)_\mathcal{P}$, $Z_\mathcal{P}(f)=Z_\mathcal{P}(|f|)=Z_\mathcal{P}(f\wedge \boldsymbol{1})=Z_\mathcal{P}(f^n)$, for all $n\in \mathbb{N}$.
		\item $Z_\mathcal{P}(\boldsymbol{0})=X$ and $Z_\mathcal{P}(\boldsymbol{1})=\emptyset$.
		\item For $f,g\in C(X)_\mathcal{P}$, $Z_\mathcal{P}(fg)=Z_\mathcal{P}(f)\cup Z_\mathcal{P}(g)$ and $Z_\mathcal{P}(|f|+|g|)=Z_\mathcal{P}(f^2+g^2)=Z_\mathcal{P}(f)\cap Z_\mathcal{P}(g)$.
		\item For $f\in C(X)_\mathcal{P}$, $r,s\in \mathbb{R}$, sets of the form $\{x\in X\colon f(x)\geq r \}$ and $\{x\in X\colon f(x)\leq s \}$ are zero sets as we have : 
		\begin{enumerate}
			\item $\{x\in X\colon f(x)\geq r \}=Z_\mathcal{P}((f-\boldsymbol{r})\wedge 0)$, and
			\item $\{x\in X\colon f(x)\geq s \}=Z_\mathcal{P}((f-\boldsymbol{s})\vee 0)$.
		\end{enumerate}
	\end{enumerate}
	
	From the above observations, it follows that $Z_\mathcal{P}[X]$ is closed under finite union and intersection.
	We know that the zero sets in $Z[X]$ (zero sets of real valued continuous functions on $X$) are $G_\delta$-sets. The following result generalises this fact for functions in $C(X)_\mathcal{P}$.
	
	\begin{theorem} \label{2.1}
		Every zero set in $X$ can be expressed as a disjoint union of a $G_\delta$-set and a set $A$ such that $\overline{A}\in \mathcal{P}$.
	\end{theorem}
	
	\begin{proof}
		Let $f\in C(X)_\mathcal{P}$. Then $Z_\mathcal{P}(f) = G\cup A$, where $G = Z_\mathcal{P}(f) \cap (X\setminus D_f)$ and $A = Z_\mathcal{P}(f) \cap D_f$ . Then $G$ is a zero set of the continuous map $f|_{X\setminus D_f}$ and is therefore a $G_\delta$-set in $X\setminus D_f$ . Since $D_f$ is an $F_\sigma$-set in $X$, $X\setminus D_f$ is a $G_\delta$-set in $X$. Thus $G$ is a $G_\delta$-set in $X$. On the other hand, $\overline{A}\subseteq \overline{D_f}\in \mathcal{P}\implies$ $\overline{A}\in \mathcal{P}$.
	\end{proof}
	
	We know that all zero sets in $Z(X)$ are closed sets in $X$. However, all zero sets in $Z_\mathcal{P}[X]$ may not be closed. In fact, we establish a necessary and sufficient condition for this to happen.
	
	\begin{theorem} \label{2.15}
		For a $\tau \mathcal{P}$-space $X$, all zero sets in $Z_\mathcal{P}[X]$ are closed if and only if $C(X)_\mathcal{P}=C(X)$.
	\end{theorem}
	
	\begin{proof}
		Let $C(X)_\mathcal{P}\neq C(X)$. Then tracing the steps in Theorem \ref{2.03}, there exists a non-isolated point $x_0\in X$ such that $\chi_{\{x_0 \}}\in C(X)_\mathcal{P}\setminus C(X)$ and $Z_\mathcal{P}(\chi_{\{x_0 \}})=X\setminus \{x_0\}$, which is not closed in $X$.
		The converse is clear.
	\end{proof}
	
	The following observations are immediate.
	\begin{proposition} \label{2.4}
		Let $X$ be a $\tau \mathcal{P}$-space. Then the following statements are true.
		\begin{enumerate}
			\item $C(X)_\mathcal{P}$ is reduced. (A commutative ring $R$ is called reduced if it does not contain any non-zero nilpotent elements.)
			\item $f$ is a unit in $C(X)_\mathcal{P}$ if and only if $Z_\mathcal{P}(f)=\emptyset$.
		\end{enumerate}
	\end{proposition}

	\begin{definition}
		Two subsets $A$ and $B$ are said to be $\mathcal{P}$-completely separated if there exists $f\in C(X)_\mathcal{P}$ such that $f(X)\subseteq [0,1]$, $f(A)=\{0\}$ and $f(B)=\{1\}$.
	\end{definition}
	Thus $\mathcal{P}$-completely separated sets reduce to completely separated sets in \cite{GJ1976}, when $\mathcal{P}=\{\emptyset \}$.

	Let $\tau_u$ be the usual topology on the set $\mathbb{R}$ of all real numbers and $\mathcal{P}=\{\emptyset, \{0\} \}$.
	Define $f\colon \mathbb{R} \rightarrow \mathbb{R}$ by: 
	$f(x)=\begin{cases}
		0 \quad &x\leq 0 \\ 1 \quad &x>0
	\end{cases}$
	Then $f\in C(X)_\mathcal{P}$. Thus $(-\infty,0]$ and $(0,\infty)$ are $\mathcal{P}$-completely separated. However $(-\infty,0]$ and $(0,\infty)$ are not completely separated.

	\begin{proposition} \label{2.2}
		Two disjoint subsets $A$ and $B$ of a $\tau \mathcal{P}$-space $X$ are $\mathcal{P}$-completely separated in $X$ if and only if they are contained in disjoint zero sets in $X$.
	\end{proposition}
	
	\begin{proof}
		The first part of the theorem can be proved by closely following the proof of Theorem 1.2 in \cite{GJ1976}.
		To prove the converse, let $A$ and $B$ be two disjoint subsets of a $\tau \mathcal{P}$-space $X$ that  are contained in disjoint zero sets, $Z_\mathcal{P}(f)$ and $Z_\mathcal{P}(g)$ respectively in $X$. Define \begin{equation*}
			h(x)=\frac{|f|(x)}{(|f|+|g|)(x)}, \quad \forall x\in X.
		\end{equation*} Then $\overline{D_{h}}\subseteq \overline{D_f}\cup \overline{D_g} \implies \overline{D_h}\in \mathcal{P}$ and so $h\in C(X)_\mathcal{P}$. Also, $h(A)=\{0\}$ and $h(B)=\{1\}$. Hence $A$ and $B$ are $\mathcal{P}$-completely separated in $X$.
	\end{proof}

	\begin{theorem} \label{2.3}
		Two disjoint subsets $A$ and $B$ of a $\tau \mathcal{P}$-space $X$ are $\mathcal{P}$-completely separated in $X$ if and only if there exists a $P\in \mathcal{P}$ such that $A\setminus P$ and $B\setminus P$ are completely separated in $X\setminus P$.
	\end{theorem}
	
	\begin{proof}
		The proof of this theorem is analogous to that of Proposition 2.3 of \cite{GGT2018}.
	\end{proof}
	
	We now introduce the concept of $z_\mathcal{P}$-filters on $X$ and $z_\mathcal{P}$-ideals of $C(X)_\mathcal{P}$.
	
	\begin{definition}
		\ \begin{enumerate}
			\item 	A non-empty subcollection $\mathcal{F}$ of $Z_\mathcal{P}[X]$ is called a $z_\mathcal{P}$-filter on $X$ if $\mathcal{F}$ satisfies the following conditions: 
			\begin{enumerate}
				\item $\emptyset \notin \mathcal{F}$,
				\item $\mathcal{F}$ is closed under finite intersections, and
				\item If $Z_1\in \mathcal{F}$ and $Z_2\in Z_\mathcal{P}[X]$ with $Z_1\subseteq Z_2$, then $Z_2\in \mathcal{F}$.
			\end{enumerate}
			\item A $z_\mathcal{P}$-filter on $X$ is said to be a $z_\mathcal{P}$-ultrafilter on $X$ if it is not properly contained in any other $z_\mathcal{P}$-filter on $X$.
			\item An ideal $I$ of $C(X)_\mathcal{P}$ is called a $z_\mathcal{P}$-ideal if $Z_\mathcal{P}^{-1}Z_\mathcal{P}[I]=I$.
		\end{enumerate}
		
	\end{definition}
	
	Using the same arguments as used in Theorems 2.3, 2.5, 2.9 and 2.11 of \cite{GJ1976}, we can prove the following theorems:
	\begin{theorem} \label{2.5}
		\	 \begin{enumerate}
			\item For any ideal $I$ of $C(X)_\mathcal{P}$, $Z_\mathcal{P}[I]$ is a $z_\mathcal{P}$-filter.
			\item For a $z_\mathcal{P}$-filter $\mathcal{F}$, $Z_\mathcal{P}^{-1}(\mathcal{F})=\{f\in C(X)_\mathcal{P}\colon Z_\mathcal{P}(f)\in \mathcal{F} \}$ is an ideal of $C(X)_\mathcal{P}$.
			\item If $M$ is a maximal ideal in $C(X)_\mathcal{P}$, then $Z_\mathcal{P}[M]$ is a $z_\mathcal{P}$-ultrafilter on $X$.
			\item For a $z_\mathcal{P}$-ultrafilter $\mathcal{U}$ on $X$, $Z_\mathcal{P}^{-1}(\mathcal{U})$ is a maximal ideal of $C(X)_\mathcal{P}$.
		\end{enumerate}
	\end{theorem}
	
	\begin{theorem} \label{2.6}
		For any $z_\mathcal{P}$-ideal $I$ of $C(X)_\mathcal{P}$, the following are equivalent:
		\begin{enumerate}
			\item $I$ is prime.
			\item $I$ contains a prime ideal.
			\item For all $g$, $h\in C(X)_\mathcal{P}$, if $gh=\boldsymbol{0}$, then $g\in I$ or $h\in I$.
			\item For every $f\in C(X)_\mathcal{P}$, there is a zero-set in $Z[I]$ on which $f$ does not change its sign.
		\end{enumerate}
	\end{theorem}
	
	\begin{theorem} \label{2.7}
		Every prime ideal in $C(X)_\mathcal{P}$ is contained in a unique maximal ideal.
	\end{theorem}
	
	\begin{corollary} \label{2.75}
		$C(X)_\mathcal{P}$ is a $pm$-ring (A commutative ring $R$ with unity is called a $pm$-ring if every prime ideal of $R$ is contained in an unique maximal ideal of $R$.).
	\end{corollary}
	
	Now we are interested to study the maximal ideal space of $C(X)_\mathcal{P}$. Let $Max(C(X)_\mathcal{P})$ be the collection of all maximal ideals of $C(X)_\mathcal{P}$.
	For $f\in C(X)_\mathcal{P}$, set
	$\mathcal{M}_f=\{M\in Max(C(X)_\mathcal{P}) \colon f\in M \}$. It is easy to see that
	$\mathcal{B}=\{\mathcal{M}_f \colon f\in C(X)_\mathcal{P} \}$ is a base for closed sets for some topology on $Max(C(X)_\mathcal{P})$, which is known as the hull-kernel topology.
	
	Let $\beta_\mathcal{P}X$ be the index set for all $z_\mathcal{P}$-ultrafilters on $X$ with the condition that for $p\in X$, $A^p=A_p=\{Z\in Z_\mathcal{P}[X]\colon p\in Z \}$. For $Z\in Z_\mathcal{P}[X]$, we set $\overline{Z}=\{p\in \beta_\mathcal{P}X \colon Z\in A^p \}$. Then $\mathcal{B'}=\{\overline{Z}\colon Z\in Z_\mathcal{P}[X] \}$ forms a base for closed sets for some topology on $\beta_\mathcal{P}X$. 
	
	The following observations  for $Z\in Z_\mathcal{P}[X]$ are immediate : 
	\begin{theorem} \label{2.8}
		\begin{enumerate}
			\item $\overline{Z}\cap X=Z$.
			\item $\overline{Z}$=$cl_{\beta_\mathcal{P}X}Z$.
		\end{enumerate}
	\end{theorem}

	\begin{theorem} \label{2.9}
		$Max(C(X)_\mathcal{P})$ is homeomorphic to $\beta_\mathcal{P}X$
	\end{theorem}
	
	\begin{proof}
		Define $\phi \colon \beta_\mathcal{P}X \longrightarrow Max(C(X)_\mathcal{P})$ by $\phi(p)=Z_\mathcal{P}^{-1}[A^p]=M^p$(say).
		Then $M^p$ is a maximal ideal of $C(X)_\mathcal{P}$. Also, $\phi$ is a bijective map, by Theorem \ref{2.5}. Let $Z=Z_\mathcal{P}(f)\in Z_\mathcal{P}[X]$. Then $\phi(\overline{Z})=\mathcal{M}_f \in \mathcal{B}$ and $\phi^{-1}(\mathcal{M}_f)=\overline{Z_\mathcal{P}(f)} \in \mathcal{B'}$. Therefore $\phi$ exchanges basic closed sets of $Max(C(X))_\mathcal{P}$ and $\beta_\mathcal{P}X$. Hence $\phi$ is a homeomorphism.
	\end{proof}
	
	Since $C(X)_\mathcal{P}$ contains unity, $Max(C(X)_\mathcal{P})$ is compact and hence $\beta_\mathcal{P}X$ is compact. Since $C(X)_\mathcal{P}$ is a $pm-$ring (by Theorem \ref{2.7}), by Theorem 1.2 of \cite{MO1971}, $\beta_\mathcal{P}X$ is Hausdorff.
	
	It can be easily seen that when $\mathcal{P}=\{\emptyset \}$, we have $\beta_\mathcal{P} X=\beta X$.
	
	We further note for $X=(0,1)\cup \{2\}$ equipped with the subspace topology inherited from the usual topology of $\mathbb{R}$ and $\mathcal{P}$ is the ideal of all closed subsets of $(0,1)$. Then, $C(X)_\mathcal{P}=\mathbb{R}^X=C(X_d)$, where $X_d=X$ equipped with the discrete topology on $X$. It is clear that $\beta_\mathcal{P}X$ has uncountably many isolated points, but $\beta X$ has only one isolated point. Thus $\beta_\mathcal{P}X$ is not homeomorphic to $\beta X$ in general.

	\section{When is $C(X)_\mathcal{P}$ closed under uniform limit?}
	
	\begin{definition}
		A sequence of functions $\{f_n \}$ in a subring $S$ of $\mathbb{R}^X$ is said to converge uniformly to a function $f$ on $X$ if for a given $\epsilon>0$, there exists $N\in \mathbb{N}$ such that $|f_n(x)-f(x)|<\epsilon$ for all $n\geq N$ and for all $x\in X$. \\ A subring $S$ of $\mathbb{R}^X$ is said to be closed under uniform limit if whenever $\{f_n\}\subseteq S$ converges uniformly to a function $f\in \mathbb{R}^X$, $f\in S$. \\
		A $\tau \mathcal{P}$-space $(X,\tau,\mathcal{P})$ is said to be a $\tau\mathcal{PU}$-space if $C(X)_\mathcal{P}$ is closed under uniform limit.
	\end{definition}
	
	It can be easily observed here that if $\mathcal{P}=\{\emptyset\}$, then $X$ is a $\tau \mathcal{PU}$-space. Another trivial example of a ring that is closed under uniform limit is $\mathbb{R}^X$. So if $C(X)_\mathcal{P}=\mathbb{R}^X$, then $X$ is a $\tau \mathcal{PU}$-space. Also, every function in $\mathbb{R}^X$ is continuous on all isolated points of $X$. In light of these observations, we have the following theorem.
	
	\begin{theorem} \label{3.1}
		Let $\mathcal{P}$ be an ideal of closed subsets of $X$ such that the set of all non-isolated points in $X$ is a member of $\mathcal{P}$. Then $C(X)_\mathcal{P}=\mathbb{R}^X$, and hence $C(X)_\mathcal{P}$ is closed under uniform limit, that is $X$ is a $\tau\mathcal{PU}$-space.
	\end{theorem}

	The converse of the above theorem holds when $\mathcal{P}=\mathcal{P}_f$, as seen in Theorem 2.9 \cite{MBM2021}.
	
	The converse of the above theorem also holds for a metrizable space $X$ and for $\mathcal{P}=\mathcal{K}$ .
	\begin{theorem} \label{3.2}
		Let $X$ be a metrizable space and $\mathcal{P}=\mathcal{K}$. If $X$ is a $\tau \mathcal{PU}$-space, then the set of all non-isolated points in $X$ is a member of $\mathcal{K}$.
	\end{theorem}
	
	\begin{proof}
		Let T be the set of all non-isolated points of $X$.
		Assume that $T$ is non-compact. Then $T$ is not sequentially compact. So, $\exists$ a sequence $\{a_n\} \in T$ which has no convergent subsequence. Set $A=\{a_n \colon n\in \mathbb{N}\}$. Then $A$ is a closed non-compact subset of $X$.
		
		For each $m\in \mathbb{N}$, define $f_m$ on $X$ as follows :
		\begin{align*}
			f_m(x)=\begin{cases}
				\frac{1}{n} & x=a_n, n<m \\
				0 & otherwise.
			\end{cases}
		\end{align*} Then, $f_m \in C(X)_K$, for each $m\in \mathbb{N}$ and $\{f_m\}$ is uniformly convergent to a function $f:X\longrightarrow \mathbb{R}$ where \begin{equation*}
			f(x)=\begin{cases}
				\frac{1}{n} & x=a_n \\ 0 & otherwise.
			\end{cases}
		\end{equation*}	Clearly, $\overline{D_f}=\overline{A}=A \notin \mathcal{K}$. Thus $f\notin C(X)_K$. Hence $C(X)_K$ is not closed under uniform limit. This completes the proof.
	\end{proof}
	
	It is well known that $C(Y)$ is closed under uniform limit for any topological space $Y$. The next natural question is that if $C(X)_\mathcal{P}$ is isomorphic to $C(Y)$, then can we conclude that $C(X)_\mathcal{P}$ is also closed under uniform limit?
	
	\begin{theorem} \label{3.3}
		Let $\mathcal{P}$ be an ideal of closed subsets of a space $X$. If $C(X)_\mathcal{P}$ is isomorphic to $C(Y)$ for some topological space $Y$, then $X$ is a $\tau \mathcal{PU}$-space.
	\end{theorem}
	
	\begin{proof}	
		Let $\phi \colon C(X)_\mathcal{P} \longrightarrow C(Y)$ be an isomorphism. First, we attempt to show that $\phi$ is an order preserving mapping.
		For that, let $g\in C(X)_\mathcal{P}$ be such that $g\geq 0$. Then $g=l^2$ for some $l\in C(X)_\mathcal{P}$. Thus $\phi(g)=\phi(l^2)=(\phi(l))^2\geq 0$.
		So, $\phi$ is order preserving.	For $f\in C(X)_\mathcal{P}$, $\phi(|f|)\geq 0$.
		Also, $(\phi(|f|))^2=\phi(|f|^2)=\phi(f^2)=(\phi(f))^2$ which implies $\phi(|f|)=|\phi(f)|$, for all $f\in C(X)_\mathcal{P}$............(1). Since $\phi$ is an isomorphism, for any rational number $r$, $\phi(\boldsymbol{r})=\boldsymbol{r}$................(2). On using the above arguments, we can show that (1) and (2) also hold for $\phi^{-1}$ as well. \\ Let $\{f_n\}$ be a sequence in $C(X)_\mathcal{P}$ converging uniformly to a function $f\in \mathbb{R}^X$. We now show that $f\in C(X)_\mathcal{P}$. 
		
		Let $\epsilon>0$ be an arbitrary rational. Then there exists $k\in \mathbb{N}$ such that $|f_n-f_m|<\boldsymbol{\epsilon}$ for all $n,m\geq k$.
		Since $\phi$ is order preserving, and using (1) and (2), we have $|\phi(f_n)-\phi(f_m)|<\phi(\boldsymbol{\epsilon})=\boldsymbol{\epsilon}$, for all $n,m\geq k$. So, there exists $h\in C(Y)$ such that $\{\phi(f_n)\}$ converges uniformly to $h\in C(Y)$. By hypothesis, $\phi$ is onto. Therefore there exists $g\in C(X)_\mathcal{P}$ such that $\phi(g)=h$.
		
		Now, for a given rational $\epsilon >0$, as $\{\phi(f_n)\}$ converges uniformly to $h\in C(Y)$, there exists $k\in \mathbb{N}$ such that \begin{align*}
			&\quad |\phi(f_n)-h|<\boldsymbol{\epsilon} \quad \forall n\geq k
			\\ &\implies |\phi(f_n)-\phi(g)|<\boldsymbol{\epsilon} \quad \forall n\geq k
			\\ &\implies |f_n-g|=|\phi^{-1}(\phi(f_n))-\phi^{-1}\phi(g)| <\phi^{-1}(\boldsymbol{\epsilon})=\boldsymbol{\epsilon} \quad \forall n\geq k.
		\end{align*}
		Thus $\{f_n\}$ converges uniformly to the function $g$. But $\{f_n\}$ converges uniformly to $f$. Thus $f=g\in C(X)_\mathcal{P}$.
		Hence $C(X)_\mathcal{P}$ is closed under uniform limit.
		
	\end{proof}
	\begin{corollary} \label{3.4}
		In a metric space $X$ the following statements are equivalent :
		\begin{enumerate}
			\item $C(X)_K$ is closed under uniform limit.
			\item The set of all non-isolated points in $X$ is compact.
			\item $C(X)_K$ is isomorphic to $C(Y)$ for some topological space $Y$.
		\end{enumerate}
	\end{corollary} 
	\begin{proof}
		(1)$\implies$(2) follows from Theorem \ref{3.2}. \\
		From Theorem \ref{3.1}, we get (2)$\implies C(X)_K=\mathbb{R}^X=C(X_d) $, where $X_d$ denotes the space $X$ with discrete topology. This proves (3). \\
		(3)$\implies$(1) follows from Theorem \ref{3.3} for $\mathcal{P}=\mathcal{K}$.
		
	\end{proof}
	
	If we consider the ideal $\mathcal{P}=\mathcal{P}_f$ on a topological space $X$, then we have the following theorem where we provide an alternative proof for Theorem 3.4 of \cite{GGT2018}.
	\begin{theorem} \label{3.5}
		For a topological space $X$, the following statements are equivalent:
		\begin{enumerate}
			\item $C(X)_F$ is closed under uniform limit.
			\item $X$ has finitely many non-isolated points.
			\item $C(X)_F$ is isomorphic to $C(Y)$, for some topological space $Y$.
		\end{enumerate}
	\end{theorem}
	\begin{proof}
		(1)$\iff$(2) follows from Theorem 2.9 in \cite{MBM2021}. \\
		From Theorem \ref{3.1}, we get (2)$\implies C(X)_F=\mathbb{R}^X=C(X_d) $, where $X_d$ denotes the space $X$ with discrete topology. This proves (3). \\
		(3)$\implies$(1) follows from Theorem \ref{3.3} for $\mathcal{P}=\mathcal{P}_f$.
	\end{proof}
	
	The fact that $C(X)$ is closed under uniform limit is used extensively to prove the Urysohn's Extension Theorem (Theorem 1.17 in \cite{GJ1976}). Our aim is to achieve an analog of that result. We need the following definitions to do this.
	
	\begin{definition}
		A subspace $S$ of $X$ is said to be $C_\mathcal{P}$-embedded if every $f\in C(S)_{\mathcal{P}_S}$ can be extended to a function in $C(X)_\mathcal{P}$. \\
		A subspace $S$ of $X$ is said to be $C_\mathcal{P}^*$-embedded if every $f\in C^*(S)_{\mathcal{P}_S}$ can be extended to a function in $C^*(X)_\mathcal{P}$.
	\end{definition}
	
	The following theorem is a generalisation of the Urysohn's Extension Theorem (Theorem 1.17 in \cite{GJ1976}) and can be proved by closely following the proof of Theorem 1.17 in \cite{GJ1976}.
	
	\begin{theorem} \label{3.6}
		Let $X$ be a $\tau \mathcal{PU}$-space. Then a subspace $S$ of $X$ is $C_\mathcal{P}^*$-embedded in $X$ if and only if any two $\mathcal{P}_S$-completely separated sets in $S$ are $\mathcal{P}$-completely separated sets in $X$.
	\end{theorem}
	
	Further, as is seen in case of $C(X)$, a $C^*_\mathcal{P}$-embedded subspace of $X$ may not be $C_\mathcal{P}$-embedded.
	
	\begin{theorem} \label{3.7}
		A $C^*_\mathcal{P}$-embedded subspace of $X$ is $C_\mathcal{P}$-embedded if and only if it is $\mathcal{P}$-completely separated from every zero set disjoint from it.
	\end{theorem}
	This can be proved by closely following the proof of Theorem 1.18 in \cite{GJ1976}.

	\section{$\tau \mathcal{P}$-compact spaces.}
	In this section, we define $\tau \mathcal{P}$-compact, $\tau \mathcal{P}$-pseudocompact and real $\tau \mathcal{P}$-compact spaces and discuss their characterisations and properties.
	\begin{definitions}
		\	\begin{enumerate}
			\item A $\tau\mathcal{P}$-space $(X,\tau,\mathcal{P})$ is said to be $\tau \mathcal{P}$-compact if for every family of zero sets $\mathcal{F}\subseteq Z_\mathcal{P}[X]$ having the finite intersection property, $\bigcap \mathcal{F}\neq \emptyset$.
			\item A $z_\mathcal{P}$-filter $\mathcal{F}$ is said to be fixed if $\bigcap \mathcal{F}\neq \emptyset$, otherwise it is said to be free. An ideal $I$ of $C(X)_\mathcal{P}$ is said to be fixed if $\bigcap Z_\mathcal{P}[I]\neq \emptyset$, otherwise it is said to be free.
			\item A $\tau \mathcal{P}$-space $X$ is said to be $\tau\mathcal{P}$-pseudocompact if $C(X)_\mathcal{P}=C^*(X)_\mathcal{P}$.
		\end{enumerate}
		
	\end{definitions}
	The following observation is immediate.
	\begin{observation}
		Let $(X,\tau,\mathcal{P})$ be a $\tau \mathcal{P}$-space such that $(X,\tau)$ is a Tychonoff space. Then $(X,\tau)$ is compact if $(X,\tau,\mathcal{P})$ is $\tau \mathcal{P}$-compact.
	\end{observation}
	However if $X=\{0 \}\cup \{\frac{1}{n}\colon n\in \mathbb{N} \}$ is endowed with the subspace topology induced from the usual topology on $\mathbb{R}$ and $\mathcal{P}=\{\emptyset, \{0 \} \}$, then $X$ is compact but is not $\tau \mathcal{P}$-compact. This shows that the converse of the above observation may fail.
	
	Pertaining to the above definitions, we have the following result that has a proof similar to that of Theorem 4.11 in \cite{GJ1976}.
	\begin{theorem} \label{4.1}
		The following statements are equivalent:
		\begin{enumerate}
			\item $X$ is $\tau \mathcal{P}$-compact.
			\item Every $z_\mathcal{P}$-filter is fixed.
			\item Every $z_\mathcal{P}$-ultrafilter is fixed.
			\item Every ideal in $C(X)_\mathcal{P}$ is fixed.
			\item Every maximal ideal in $C(X)_\mathcal{P}$ is fixed.
		\end{enumerate}
	\end{theorem}

	\begin{theorem} \label{4.2}
		If $(X,\tau,\mathcal{P})$ is $\tau \mathcal{P}$-compact, then it is $\tau \mathcal{P}$-pseudocompact.
	\end{theorem}
	
	\begin{proof}
		Let $(X,\tau,\mathcal{P})$ be $\tau\mathcal{P}$-compact.
		If possible, let there exists an unbounded function $f\in C(X)_\mathcal{P}$ and consider $Z_n=f^{-1}((-\infty,-n]\cup [n,\infty))$, for each $n\in \mathbb{N}$. Then $\mathcal{B}=\{Z_n \colon n\in \mathbb{N} \}$ is a subcollection of $Z_\mathcal{P}[X]$. Each $Z_n$ is non-empty as $f$ is unbounded. Thus $\exists$ a $z_\mathcal{P}$-filter $\mathcal{F}$ on $X$ such that $\mathcal{B}\subseteq \mathcal{F}$. By Theorem \ref{4.1}, $\emptyset \neq \bigcap \mathcal{F}\subseteq \bigcap \mathcal{B}$. Let $x\in \bigcap \mathcal{B}$. Then $|f(x)|\geq n$, for all $n\in \mathbb{N}$, which is not possible.
		Thus $C^*(X)_\mathcal{P}=C(X)_\mathcal{P}$. Hence $X$ is $\tau \mathcal{P}$-pseudocompact.
	\end{proof}
	
	It is important to note that the converse of the above result is false in general, as is shown in the following example.
	
	\begin{counter example} \label{4.3}
		Let us consider $X=\mathbb{R}$ equipped with the co-countable topology $\tau$ and take $\mathcal{P}=\mathcal{P}_f$, the collection of all finite subsets of $X$.	For any $f\in C(X)_\mathcal{P}$, we have $f$ is continuous on $X\setminus D_f$, which is open in $X$, as $D_f$ is finite. So $f$ is constant on $X\setminus D_f$. Further, $D_f$ is finite. Thus $f\in C^*(X)_\mathcal{P}$. Therefore $(X,\tau,\mathcal{P})$ is $\tau \mathcal{P}$-pseudocompact.
		For each $x\in X$, let us consider $f_x=\chi_{_{X\setminus \{x\}} }\in C(X)_\mathcal{P}$. Then $Z_\mathcal{P}(f_x)=\mathbb{R}\setminus \{x\}$. Let $\mathcal{B}=\{Z_\mathcal{P}(f_x) \colon x\in X \}\subseteq Z_\mathcal{P}[X]$. Then $\mathcal{B}$ has the finite intersection property. However $\bigcap \mathcal{B}=\emptyset$. Thus $X$ is not $\tau \mathcal{P}$-compact.

	\end{counter example}

	An interesting characterisation of a $\tau \mathcal{PU}$-space can be achieved for a $\tau \mathcal{P}$-pseudocompact space. It has already been noted in \cite{GJ1976} that $C^*(X)$, equipped with the norm $||f||=\sup \{|f(x)|\colon x\in X \}$ is a Banach space. However, $C^*(X)_\mathcal{P}$ equipped with the same norm may not be a Banach space.
	
	Since the uniform convergence of a sequence of functions in $C^*(X)_\mathcal{P}$ is same as the convergence of the sequence under the norm defined above, we have the following result.
	
	\begin{theorem} \label{3.0}
		A $\tau\mathcal{P}$-pseudocompact space $X$ is a $\tau \mathcal{PU}$-space if and only if $(C^*(X)_\mathcal{P},||.||)$ is a Banach space.
	\end{theorem}

	In order to establish a concept similar to real compactness of a topological space $X$, we need the following definitions.
	\begin{definitions}
		\begin{enumerate}
			\item An ideal $I$ of a partially ordered ring $R$ is said to be convex if  $0\leq x \leq y$ with $y\in I$ implies $x\in I$.
			\item An ideal $I$ of a lattice ordered ring $R$ is said to be absolutely convex if whenever $|x|\leq |y|$ with $y\in I$, then $x\in I$.
			\item A totally ordered field $R$ is said to be archimedean if for every element $a\in R$, there exists $n\in \mathbb{N}$ such that $\boldsymbol{n}\geq a$, where $\boldsymbol{n}$ denotes $1_R$(identity of $R$) added $n$ times. 
		\end{enumerate}
	\end{definitions}
	
	\begin{proposition} \label{4.4}
		Every $z_\mathcal{P}$-ideal of the ring $C(X)_\mathcal{P}$ is absolutely convex.
	\end{proposition}
	
	\begin{proof}
		Straightforward.
	\end{proof}
	
	From Theorems 5.2, 5.3 in \cite{GJ1976}, we have the following theorem.
	\begin{theorem} \label{4.5}
		Let $I$ be an absolutely convex ideal in a lattice ordered ring $A$. Then $A/I$ is a lattice ordered ring, according to the following definition: $I(a)\geq 0$ if there exists $x\in A$ such that $x\geq 0$ and $I(a)=I(x)$.
	\end{theorem}
	The following proposition is immediate and can be proved by closely following the discussions in 5.4(a) and 5.4(b) of \cite{GJ1976}.
	\begin{proposition} \label{4.6}
		For a $z_\mathcal{P}$-ideal $I$, \begin{enumerate}
			\item $I(f)\geq 0$ if and only if there exists $Z\in Z_\mathcal{P}[I]$ such that $f(x)\geq 0$ for all $x\in Z$. 
			\item If $f>\boldsymbol{0}$ on some zero set of a function in $I$ then $I(f)>0$.
		\end{enumerate}
	\end{proposition}

	\begin{definition}
		A maximal ideal $M$ of the ring $C(X)_\mathcal{P}$ is said to be real if $C(X)_\mathcal{P}/ M$ is isomorphic to the field of real numbers.
	\end{definition}

	\begin{proposition} \label{4.7}
		A maximal ideal $M$ of $C(X)_\mathcal{P}$ is real if and only if $C(X)_\mathcal{P} /M$ is archimedean.
	\end{proposition}
	
	\begin{proof}
		Let $C(X)_\mathcal{P}/M$ be archimedean. Then by Theorem 0.21 in \cite{GJ1976}, there exists an isomorphism $\phi$ from $C(X)_\mathcal{P}/M$ to a subfield $S$ of $\mathbb{R}$. Also, $\psi \colon \mathbb{R} \longrightarrow C(X)_\mathcal{P}/M$ defined by $\psi(r)=M(\boldsymbol{r})$ for all $r\in \mathbb{R}$ is an injective homomorphism. Thus $\phi \circ \psi$ is an injective homomorphism from $\mathbb{R}$ onto a subfield of $\mathbb{R}$. By Theorem 0.22 in \cite{GJ1976}, $\phi \circ \psi$ is the identity homomorphism from $\mathbb{R}$ onto $\mathbb{R}$. Thus $\phi (C(X)_\mathcal{P}/M)=\mathbb{R}\implies$ $M$ is real. The converse is obvious.
	\end{proof}
	
	\begin{definition}
		A $\tau\mathcal{P}$-space $(X,\tau,\mathcal{P})$ is called $\tau\mathcal{P}$-real compact if every real maximal ideal of $C(X)_\mathcal{P}$ is fixed.
	\end{definition}
	
	Real maximal ideals of $C(X)_\mathcal{P}$ has the following characterization for a $\tau \mathcal{PU}$-spaces :
	\begin{theorem} \label{4.8}
		Let $X$ be a $\tau \mathcal{PU}$-space. Then the following statements are equivalent.
		\begin{enumerate}
			\item $M$ is a real maximal ideal of $C(X)_\mathcal{P}$.
			\item $Z_\mathcal{P}[M]$ is closed under countable intersection.
			\item $Z_\mathcal{P}[M]$ has the countable intersection property.
		\end{enumerate}
	\end{theorem}
	The above result can be proved by following the proof of Theorem 5.14 \cite{GJ1976}.
	
	Next we state a characterization of $\tau \mathcal{P}$-compact spaces using real $\tau \mathcal{P}$-compact and $\tau \mathcal{P}$-pseudocompact spaces.
	
	\begin{theorem} \label{4.9}
		A $\tau\mathcal{P}$-space $X$ is $\tau\mathcal{P}$-compact if and only if it is both $\tau\mathcal{P}$-pseudocompact and $\tau\mathcal{P}$-real compact.
	\end{theorem}
	
	\begin{proof}
		Let $X$ be both $\tau\mathcal{P}$-pseudocompact and $\tau\mathcal{P}$-real compact. Let $\mathcal{U}$ be a $z_\mathcal{P}$-ultrafilter on $X$ and $M=Z_\mathcal{P}^{-1}(\mathcal{U})$. Then $M$ is a maximal ideal of $C(X)_\mathcal{P}$.	If $M$ is not a real maximal ideal, then by Proposition \ref{4.7}, there exists $f\in C(X)_\mathcal{P}$ such that $M(f)\geq M(\boldsymbol{n})$ for all $n\in \mathbb{N}$. So, for each $n\in \mathbb{N}$, there exists $x_n\in X$ such that $f(x_n)\geq n$.	This contradicts that $X$ is $\tau\mathcal{P}$-pseudocompact. Thus $M$ is a real maximal ideal. Since $X$ is $\tau\mathcal{P}$-real compact, $\bigcap Z_\mathcal{P}[M]\neq \emptyset$ and so $\bigcap \mathcal{U}\neq \emptyset$. Thus $X$ is $\tau\mathcal{P}$-compact. \\ The converse follows partially from Theorem \ref{4.2} and the rest is straightforward.
		
	\end{proof}
	
	It is to note that every $\tau \mathcal{P}$-compact space is always real $\tau \mathcal{P}$-compact, the converse of this statement may not be true as is seen in the following example.
	
	\begin{counter example} \label{4.10}
		Let us consider $X=\mathbb{R}$ equipped with usual topology and $\mathcal{P}$ be any ideal of closed subsets of $X$. Then $(X,\tau,\mathcal{P})$ is a $\tau \mathcal{P}$ space. Let $M$ be a real maximal ideal of $C(X)_\mathcal{P}$. Then $\phi \colon \mathbb{R} \longrightarrow C(X)_\mathcal{P}/M$ defined as $\phi(r)=M(\boldsymbol{r})$ is an isomorphism. Now $i$ (identity map on $X$) $\in C(X)\subseteq C(X)_\mathcal{P}\implies$ $\exists$ $r\in \mathbb{R} $ such that $M(\boldsymbol{r})=M(i)$. This implies that $\{r\}\subseteq Z_\mathcal{P}[M]$, that is, the ideal $M$ is fixed. Thus $X$ is real $\tau \mathcal{P}$-compact.
		However, since $X$ is not compact, it is not $\tau \mathcal{P}$-compact.
	\end{counter example}
	
	Two interesting observations can be made here:
	\begin{enumerate}
		\item A $\tau \mathcal{P}$-pseudocompact space may not be real $\tau \mathcal{P}$-compact.
		This can be seen using Counter Example \ref{4.3} and Theorem \ref{4.9}.
		\item A real $\tau \mathcal{P}$-compact space may not be $\tau \mathcal{P}$-pseudocompact.
		This can be seen using Counter Example \ref{4.10} and Theorem \ref{4.9}.
	\end{enumerate}

	\section{When $\mathcal{P}$ contains all singleton subsets of $X$}
	
	In this section, we study properties of $C(X)_\mathcal{P}$ under the restriction that $\mathcal{P}$ contain all singleton subsets of $X$. Since $X$ is $T_1$ and all singleton subsets of $X$ are in $\mathcal{P}$, we have $\chi_{\{x\}} \in C(X)_\mathcal{P}$ for every $x\in X$. Two of these type of rings, viz, $T'(X)$ and $C(X)_F$ have been studied in \cite{GGT2018}, \cite{AK2021}, \cite{MBM2021} and \cite{A2010}. \\ Throughout this section, we assume that any ideal $\mathcal{P}$ of closed subsets of a $\tau \mathcal{P}$-space $X$ contains all singleton subsets of $X$ (unless otherwise specified).
	
	For any non-unit element $f\in C(X)_\mathcal{P}$, there exists $x_0\in X$ such that $f(x_0)=0$. This gives $f\chi_{_{X\setminus \{x_0\}}}=f$. This shows that $C(X)_\mathcal{P}$ is almost regular which is summarised in the following theorem.
	
	\begin{theorem} \label{5.01}
		For a $\tau \mathcal{P}$-space $X$, $C(X)_\mathcal{P}$ is an almost regular ring.
	\end{theorem}
	Also, for a non-unit element $f\in C(X)_\mathcal{P}$, there exists $y\in Z_\mathcal{P}(f)$ such that $f\chi_{\{y\}}=\boldsymbol{0}$. So $f$ is a zero divisor. Thus we have the following result.
	
	\begin{proposition} \label{5.1}
		For a $\tau \mathcal{P}$-space $X$, $f\in C(X)_\mathcal{P}$ is either a zero divisor or an unit.
	\end{proposition}
	The above result might fail if we remove the condition that $\mathcal{P}$ contains all singleton subsets of $X$.
	
	\begin{counter example} \label{5.2}
		Let us consider $X=\mathbb{R}$ with usual topology and $\mathcal{P}=$ the collection of all closed subsets of $(0,\infty)$. Then $f(x)=|x|$ is such that $f\in C(X)_\mathcal{P}$.
		Since $Z_\mathcal{P}(f)\neq \emptyset$, $f$ is a non-unit element. 
		Now, let $g\in C(X)_\mathcal{P}$ be such that $fg=\boldsymbol{0}$. Then as $f(x)\neq 0$ for all $x\neq 0$, we must have, $g(x)=0$ for all $x\neq 0$. \\ If $g(0)\neq 0$, then $\overline{D_g}=\{0\}\notin \mathcal{P}$, which contradicts that $g\in C(X)_\mathcal{P}$. So, $g(0)=0$ and thus $g=\boldsymbol{0}$. Therefore $f$ is not a zero divisor.
	\end{counter example}
	
	Our next aim is to generalise Proposition 3.1 of \cite{GGT2018} in the following way.
	\begin{proposition} \label{5.3}
		The following statements are equivalent for a $\tau \mathcal{P}$-space, $(X,\tau,\mathcal{P})$.
		\begin{enumerate}
			\item \label{5.31} $C(X)=C(X)_\mathcal{P}$.
			\item \label{5.32} $X$ is discrete.
			\item \label{5.33} $C(X)_\mathcal{P}$ is a ring of quotients of $C(X)$.
		\end{enumerate} 
	\end{proposition}
	\begin{proof}
		If $X$ is discrete, then all functions in $\mathbb{R}^X$ are continuous. So, $\mathbb{R}^X=C(X)=C(X)_\mathcal{P}$. Further, when $C(X)=C(X)_\mathcal{P}$, then for every $x\in X$, $\chi_{\{x\}}\in C(X)_\mathcal{P}=C(X)$ which implies that $\{x\}$ is open. Thus $X$ is discrete. This shows that \ref{5.31} and \ref{5.32} are equivalent. To show that \ref{5.33} and \ref{5.32} are equivalent, it is already seen that if $X$ is discrete, then $C(X)=C(X)_\mathcal{P}$. So, for every $f\in C(X)_\mathcal{P}$, we have $f\cdot \boldsymbol{1}=f\in C(X)_\mathcal{P}=C(X)$. This proves \ref{5.33}. Finally, let \ref{5.33} hold. Then for each $x\in X$, there exists $f\in C(X)$ such that $f\chi_{\{x \}}\in C(X)\setminus \{\boldsymbol{0} \}$. This implies that $f(x)\neq 0$. But $\chi_{\{x \}}=\frac{1}{f(x)}f\chi_{\{x \}}\in C(X)$. This shows that $X$ is discrete.
	\end{proof}
	
	Similarly, following the proof of Theorem 3.2 of \cite{GGT2018} we have:
	\begin{theorem} \label{5.4}
		The following are equivalent for a $\tau \mathcal{P}$-space, $X$:
		\begin{enumerate}
			\item $X$ is finite.
			\item Each proper ideal of $C(X)_\mathcal{P}$ is fixed.
			\item Each maximal ideal of $C(X)_\mathcal{P}$ is fixed.
			\item Each proper ideal of $C^*(X)_\mathcal{P}$ is fixed.
			\item Each maximal ideal of $C^*(X)_\mathcal{P}$ is fixed.
		\end{enumerate}
	\end{theorem}

	\begin{corollary} \label{5.5}
		A $\tau \mathcal{P}$-space $X$ is $\tau \mathcal{P}$-compact if and only if $X$ is finite.
	\end{corollary}
	\begin{proof}
		This follows from the above theorem and Theorem \ref{4.1}.
	\end{proof}
	
	However, this result may fail even if $\mathcal{P}$ fails to contain all singleton subsets of $X$.
	
	\begin{counter example} \label{5.6}
		Let $X=\{\frac{1}{n}\colon n\in \mathbb{N} \}\cup \{0\}$ be the subspace of real line and $\mathcal{P}=\{\emptyset , \{1\} \}$. Then as $1$ is an isolated point,\\ $C(X)_\mathcal{P}=C(X)$.
		Also $X$ is compact, which implies that every maximal ideal of $C(X)=C(X)_\mathcal{P}$ is fixed; even though $X$ is an infinite set. 
	\end{counter example}
	
	We next move on to discuss certain rings properties of $C(X)_\mathcal{P}$. We start by discussing the structure of the minimal ideals and the socle of the ring.
	
	\begin{theorem} \label{5.7}
		The following assertions are true for a $\tau \mathcal{P}$-space $X$.
		\begin{enumerate}
			\item \label{5.71} A non zero ideal $I$ of $C(X)_\mathcal{P}$ is minimal if and only if there exists an $\alpha \in X$ such that $I=<\chi_{\{\alpha \}}>$ if and only if $|Z_\mathcal{P}[I]|=2$.
			\item \label{5.72} The socle of $C(X)_\mathcal{P}$ consists of all functions that vanish everywhere except on a finite set.
			\item \label{5.73} The socle of $C(X)_\mathcal{P}$ is essential and free.
		\end{enumerate}
	\end{theorem}
	
	\begin{proof}
		\begin{enumerate}
			\item Let $I$ be a non zero minimal ideal of $C(X)_\mathcal{P}$. For $f\in I\setminus \{\boldsymbol{0} \}$, there exists $\alpha \in X$ such that $f(\alpha)\neq 0$. Therefore $\chi_{\{\alpha \}}=\frac{1}{f(\alpha)}\chi_{\{\alpha \}}f\in I$. Since $I$ is a minimal ideal of $C(X)_\mathcal{P}$, it follows that $I=<\chi_{\{\alpha \}}>$. This shows that $Z_\mathcal{P}[I]=\{Z_\mathcal{P}(\boldsymbol{0}), Z_\mathcal{P}(\chi_{\{\alpha \}}) \}=\{X, X\setminus \{\alpha\} \}$. Thus, $|Z_\mathcal{P}[I]|=2$. 
			
			Next we show that $<\chi_{\{\alpha \}}>$ is a minimal ideal of $C(X)_\mathcal{P}$. Let $I$ be an ideal of $C(X)_\mathcal{P}$, $\{\boldsymbol{0} \}\subsetneq I\subseteq <\chi_{\{\alpha \}}>$. Then there exists $f\in I\setminus \{\boldsymbol{0} \}\subseteq <\chi_{\{\alpha \}}>$. So $f=g\chi_{\{\alpha \}}$, for some $g\in C(X)_\mathcal{P}$. But $f=g\chi_{\{\alpha \}}=g(\alpha)\chi_{\{\alpha \}} \implies \chi_{\{\alpha \}}=\frac{1}{g(\alpha)}f\in I$. Therefore $I=<\chi_{\{\alpha \}}>$.
			
			Finally we assume that $|Z_\mathcal{P}[I]|=2$ and show that $I$ is a minimal ideal of $C(X)_\mathcal{P}$. There exists $f\in I$ such that $f(\alpha)\neq 0$ for some $\alpha \in X$. So $\chi_{\{\alpha \}}=\frac{1}{f(\alpha)}\chi_{\{\alpha \}}f\in I$. By our assumption, for any non zero function $g\in I$, $Z_\mathcal{P}(g)=Z_\mathcal{P}(\chi_{\{\alpha \}})=X\setminus \{\alpha\}$. So every non zero $g\in I$ is of the form $g=g(\alpha)\chi_{\{\alpha \}}$. Therefore $I=\{c\chi_{\{\alpha \}} \colon c\in \mathbb{R} \}=<\chi_{\{\alpha \}}>$, which is a minimal ideal, as seen above.
			
			\item By \ref{5.71}, the socle of $C(X)_\mathcal{P}$, \[Soc(C(X)_\mathcal{P})=\sum_{\alpha \in X}<\chi_{\{\alpha \}}>=< \bigl\{\chi_{\{\alpha \}}\colon \alpha \in X \bigr\} >. \] 
			Thus every function in $Soc(C(X)_\mathcal{P})$ vanishes everywhere except for a finitely many points of $X$. Conversely, let $f\in C(X)_\mathcal{P}$ be such that it vanishes everywhere except for a finitely many points, that is $Z_\mathcal{P}(f)=X\setminus \{\alpha_i \colon \alpha_i\in X, i=1,...n \}$ where $n\in \mathbb{N}$. Then \[f=\sum_{i=1}^{n}f(\alpha_i)\chi_{\{\alpha_i \}}\in Soc(C(X)_\mathcal{P}) .\]
			\item Let $I$ be a non zero ideal of $C(X)_\mathcal{P}$. Then there exists $f\in I$ such that $f(\alpha)\neq 0$ for some $\alpha \in X$. From \ref{5.71}, we have $\chi_{\{\alpha \}}\in Soc(C(X)_\mathcal{P})$ and $\chi_{\{\alpha \}}=\frac{1}{f(\alpha)}\chi_{\{\alpha \}}f\in I$. This ensures that $Soc(C(X)_\mathcal{P})\cap I\neq \emptyset$. Thus $Soc(C(X)_\mathcal{P})$ is an essential ideal. Also, for an arbitrary $\alpha \in X$, $\chi_{\{\alpha \}}\in Soc(C(X)_\mathcal{P})$ and $\chi_{\{\alpha \}}(\alpha)=1$. So $\alpha \notin Z_\mathcal{P}[Soc(C(X)_\mathcal{P})]$. This ensures that $Soc(C(X)_\mathcal{P})$ is a free ideal.

		\end{enumerate}
	\end{proof}
	
	We shall note here that the condition that $\mathcal{P}$ contains all singleton subsets of $X$ is not a necessary condition for \ref{5.72} in the above theorem. This can be seen by taking a Tychonoff space $X$ and $\mathcal{P}=\{\emptyset\}$. Here $C(X)_\mathcal{P}=C(X)$. The rest follows from Proposition 2.2 in \cite{A1997}.
	
	Using the above results, we establish a condition under which $C(X)_\mathcal{P}$ is an Artinian Ring. We need the following result to do this.
	
	\begin{proposition} \label{5.8}
		$Soc(C(X)_\mathcal{P})=C(X)_\mathcal{P}$ if and only if $X$ is finite.
	\end{proposition}
	\begin{proof}
		Let $X=\{x_1,x_2,...,x_n \}.$ Then $\boldsymbol{1}=\sum_{i=1}^{n}\chi_{\{x_i\}}\in Soc(C(X)_\mathcal{P}).$ Thus $C(X)_\mathcal{P}=Soc(C(X)_\mathcal{P})$. Conversely, let $Soc(C(X)_\mathcal{P})=C(X)_\mathcal{P}$. Then there exists $f_i\in C(X)_\mathcal{P}$ for $i=1,2,...,n$ such that $\boldsymbol{1}=\sum_{i=1}^{n}f_i\chi_{\{x_i\}}=\sum_{i=1}^{n}f_i(x_i)\chi_{\{x_i\}}.$ So, for each $x\in X$, \[1=\sum_{i=1}^{n}f_i(x_i)\chi_{\{x_i\}}(x)=f_i(x_i)\chi_{\{x_i\}}(x) \text{ for some }i\in \{1,2,...,n \}  \] which implies $\chi_{\{x_i\}}(x)=1 \implies x=x_i$. Thus $X$ is finite.
	\end{proof}
	
	\cite{BK2000} tells us that a commutative ring $R$ with unity is semisimple if and only if $rad(R)=\{0\}$. Further $R$ is Artinian semisimple if and only if $R$ equals the sum of its minimal ideals.
	
	In the ring $C(x)_\mathcal{P}$, it is easy to see that $\bigcap_{p\in X} M_p =\{\boldsymbol{0} \}$. So $rad(C(X)_\mathcal{P})=\{\boldsymbol{0} \}$. Thus $C(X)_\mathcal{P}$ is semisimple. 
	
	Under the assumption that $\mathcal{P}$ contains all singleton subsets of $X$ and using the above discussions, we have the following corollary.
	
	\begin{corollary} \label{5.9}
		$C(X)_\mathcal{P}$ is an Artinian ring if and only if $X$ is finite.
	\end{corollary}
	
	An obvious question arises here : When is $C(X)_\mathcal{P}$ Noetherian? We are able to answer this question when $\mathcal{P}$ contains all singleton subsets of $X$. 
	
	\begin{theorem} \label{5.95}
		$C(X)_\mathcal{P}$ is an Noetherian ring if and only if $X$ is finite.
	\end{theorem}
	
	\begin{proof}
		By \ref{5.9} and the fact that an Artinian commutative ring is Noetherian, $X$ is finite implies that $C(X)_\mathcal{P}$ is Noetherian. Conversely let $X$ be an infinite set. Then $X$ contains a countably infinite set $\{x_n\colon n\in \mathbb{N} \}$. Then $<\chi_{\{x_1\}}> \subsetneqq <\chi_{\{x_1,x_2\}}> \subsetneqq <\chi_{\{x_1,x_2,x_3\}}> \subsetneqq ... $ gives an unbounded ascending chain of ideals of $C(X)_\mathcal{P}$. Thus $C(X)_\mathcal{P}$ is not Noetherian.
	\end{proof}
	
	It is important to note that the condition $\mathcal{P}$ contains all singleton subsets of $X$ is not superfluous in \ref{5.7}(\ref{5.71}), \ref{5.9} and \ref{5.95}. We see that in the following example.
	
	\begin{counter example}
		Let $X=\mathbb{R}$ with cofinite topology and $\mathcal{P}=\{\emptyset \}$. Then $C(X)_\mathcal{P}=C(X)$ which consists of only the constant functions on $\mathbb{R}$. Thus $C(X)_\mathcal{P}$ is isomorphic to $\mathbb{R}$ and the only ideals of $C(X)_\mathcal{P}$ are $\{\boldsymbol{0} \}$ and itself. So $\{\boldsymbol{0} \}$ is the only minimal ideal of $C(X)_\mathcal{P}$ and is not generated by $\chi_{\{x\}}$ for any $x\in X$. Also $|Z_\mathcal{P}[\{\boldsymbol{0} \}]|=1$. Further $C(X)_\mathcal{P}$ is both Artinian and Noetherian, even though $X$ is an infinite set. 
	\end{counter example}
	We continue the study of ring properties of $C(X)_\mathcal{P}$ and establish a set of equivalent conditions to determine when is $C(X)_\mathcal{P}$ an $IN$-ring, $SA$-ring and/or a Baer ring.
	
	\begin{theorem} \label{5.10}
		The following statements are equivalent for a $\tau \mathcal{P}$-space $(X,\tau, \mathcal{P})$.
		\begin{enumerate}
			\item \label{5.101} Any two disjoint subsets of $X$ are $\mathcal{P}$-completely separated.
			\item \label{5.102} $C(X)_\mathcal{P}$ is an $IN$-ring.
			\item \label{5.103} $C(X)_\mathcal{P}$ is an $SA$-ring.
			\item \label{5.104} $C(X)_\mathcal{P}$ is an Baer ring.
			\item \label{5.105} The space of all prime ideals of $C(X)_\mathcal{P}$ is extremally disconnected.
			\item \label{5.106} Any subset of $X$ is of the form $coz(e)$ for some idempotent $e\in C(X)_\mathcal{P}$.
			\item \label{5.107} For any subset $A$ of $X$, there exists $P\in \mathcal{P}$ such that $A\setminus P$ is a clopen subset of $X\setminus P$.
		\end{enumerate}
	\end{theorem}
	
	To prove this result, we need the following two lemmas.
	
	\begin{lemma} \label{5.11}
		For any subset $A$ of $X$, there exists a subset $S$ of $C(X)_\mathcal{P}$ such that \[A=\bigcup coz[S]=\bigcup\{coz(f)\colon f\in S \}. \]
	\end{lemma}
	This follows directly since \[A=\bigcup coz[S] \text{ where } S=\{\chi_{\{x \}}\colon x\in A \} \text{ and }\chi_{\{x \}}\in C(X)_\mathcal{P} \text{ for all }x\in X. \]
	
	\begin{lemma} \label{5.12}
\		\begin{enumerate}
			\item Let $I$ and $J$ be ideals of $C(X)_\mathcal{P}$. Then $Ann(I)\subseteq Ann(J)$ if and only if $\bigcap Z_\mathcal{P}[I]\subseteq \bigcap Z_\mathcal{P}[J]$ if and only if $\bigcap coz[J]\subseteq \bigcap coz[I]$.
			\item For any subset $S$ of $C(X)_\mathcal{P}$, $Ann(S)=\{f\in C(X)_\mathcal{P}\colon \bigcup coz[S]\subseteq Z_\mathcal{P}(f) \}$.
		\end{enumerate}
	\end{lemma}
	\begin{proof}
		\begin{enumerate}
			\item Let $Ann(I)\subseteq Ann(J)$ and $x\in \bigcap Z_\mathcal{P}[I]$. Then $f(x)=0$ for all $f\in I$ $\implies \chi_{\{x \}}f=\boldsymbol{0}$ for all $f\in I$. Therefore $\chi_{\{x \}}\in Ann(I)\subseteq Ann(J)$ $\implies \chi_{\{x \}}g=\boldsymbol{0}$  for all $g\in J$. So $g(x)=0$ for all $g\in J$ $\implies x\in \bigcap Z_\mathcal{P}[J]$. Conversely, let $\bigcap Z_\mathcal{P}[I]\subseteq \bigcap Z_\mathcal{P}[J]$ and $f\in Ann(I)$. Then $fh=\boldsymbol{0}$ for all $h\in I$. So $coz(f)\subseteq \bigcap Z_\mathcal{P}[I]\subseteq \bigcap Z_\mathcal{P}[J]$. Thus $ fh_1=\boldsymbol{0}$ for all $h_1\in J$ and hence $ f\in Ann(J)$.
			\item Let $f\in Ann(S)$. Then $fg=\boldsymbol{0}$ for all $g\in S$. Therefore for $x\in \bigcup coz[S]$, $f(x)=0$. Conversely, let $f\in C(X)_\mathcal{P}$ be such that $\bigcup coz[S]\subseteq Z_\mathcal{P}(f)$ and $g\in S$. Then $coz(g)\subseteq \bigcup coz[S]\subseteq Z_\mathcal{P}(f)$. Therefore $fg=\boldsymbol{0}$. Hence $f\in Ann(S)$.
		\end{enumerate}
	\end{proof}
	We now prove Theorem \ref{5.10}.
	\begin{proof}
		Since $C(X)_\mathcal{P}$ is a reduced commutative ring, it follows from Lemma \ref{1.2} that the statements from (\ref{5.102}) to (\ref{5.105}) are equivalent.
		We use Lemma \ref{1.1} to prove (\ref{5.101}) is equivalent to (\ref{5.102}). Let (\ref{5.101}) hold and let $I$ and $J$ be orthogonal ideals of $C(X)_\mathcal{P}$. Then $\bigcup coz[I]$ and $\bigcup coz[J]$ are disjoint subsets of $X$. By (\ref{5.101}) there exists disjoint zero sets in $C(X)_\mathcal{P}$, $Z_\mathcal{P}(f_1)$ and $Z_\mathcal{P}(f_2)$ such that $\bigcup coz[I]\subseteq Z_\mathcal{P}(f_1)$ and $\bigcup coz[J]\subseteq Z_\mathcal{P}(f_2)$. This implies that $f_1\in Ann(I)$ and $f_2\in Ann(J)$. So ${f_1}^2+{f_2}^2$ is a unit in $Ann(I)+Ann(J)$. This proves (\ref{5.102}).          
		Next let (\ref{5.102}) be true and also let $A$ and $B$ be disjoint subsets of $X$. By Lemma \ref{5.11}, there exists subsets $S_A, S_A\subseteq C(X)_\mathcal{P}$ such that $A=\bigcup coz[S_A]$ and $B=\bigcup coz[S_B]$. Let $I=<S_A>$ and $J=<S_B>$. Then $\bigcup coz[I]$ and $\bigcup coz[J]$ are disjoint sets (as $A$ and $B$ are disjoint). Therefore $I$ and $J$ are orthogonal ideals of $C(X)_\mathcal{P}$. By ( \ref{5.102}) and Lemma \ref{1.1}, $Ann(I)+Ann(J)=C(X)_\mathcal{P}$. So there exists $h_1\in Ann(I)$ and $h_2\in Ann(J)$ such that $h_1+h_2=\boldsymbol{1}$, which is a unit. Therefore $Z_\mathcal{P}(h_1)$ and $Z_\mathcal{P}(h_1)$ are disjoint. Further $A=\bigcup coz[S_A]\subseteq \bigcup coz[I]\subseteq Z_\mathcal{P}(h_1)$ (since $h_1\in Ann(I)$). Similarly $B\subseteq Z_\mathcal{P}(h_2)$. This proves (\ref{5.101}).
		
		We next show that (\ref{5.104}) is equivalent to (\ref{5.106}). Let $A\subseteq X$. Then there exists $S\subseteq C(X)_\mathcal{P}$ (by Lemma \ref{5.11}) such that $A=\bigcup coz[S]$. Define $I$ to be the ideal generated by $S$. By (\ref{5.104}) there exists an idempotent $e'\in C(X)_\mathcal{P}$ such that $Ann(I)=<e'>=Ann(<e>)$, where $e=\boldsymbol{1}-e'$ is also an idempotent. By Lemma \ref{5.12}, we have $\bigcup coz[I]=\bigcup coz[<e>]$. It can be easily seen that $\bigcup coz[S]=\bigcup coz[I]$. Thus $A=\bigcup coz[S]=\bigcup coz[<e>]=X\setminus Z_\mathcal{P}(e)$. This proves (\ref{5.106}). Let (\ref{5.106}) be true and $I$ be an ideal of $C(X)_\mathcal{P}$. By (\ref{5.106}) there exists an idempotent $e\in C(X)_\mathcal{P}$ such that $\bigcup coz[I]=coz(e)$. By Lemma \ref{5.12},
		$Ann(I)=\{f\in C(X)_\mathcal{P} \colon \bigcup coz[I]\subseteq Z_\mathcal{P}(f) \}=\{f\in C(X)_\mathcal{P} \colon coz(e)\subseteq Z_\mathcal{P}(f) \} =Ann(e)=<(\boldsymbol{1}-e)>.$ This shows that $C(X)_\mathcal{P}$ is a Baer ring.
		
		Finally we show that (\ref{5.106}) and (\ref{5.107}) are equivalent Let $A\subseteq X$. By (\ref{5.106}), $A=coz(e)$ for some idempotent $e\in C(X)_\mathcal{P}$. Let $P=\overline{D_e}\in \mathcal{P}$. It is easy to see that $coz(e)=Z_\mathcal{P}(\boldsymbol{1}-e)$. Thus $A\setminus P=X\setminus Z(e|_{X\setminus P})=Z((\boldsymbol{1}-e)|_{X\setminus P})$ is clopen in $X\setminus P$. Let (\ref{5.107}) hold and $A\subseteq X$. Then by (\ref{5.107}), there exists $P\in \mathcal{P}$ such that $A\setminus P$ is clopen in $X\setminus P$. Define $e=\chi_A$. Then $e|_{A\setminus P}$ is continuous on $X\setminus P$. Therefore $\overline{D_e}\subseteq P\in \mathcal{P}$. So $e\in C(X)_\mathcal{P}$ and $A=coz(e)$.
	\end{proof}

	\section{When is $C(X)_\mathcal{P}$ regular?}
	
	We study the regularity of the ring $C(X)_\mathcal{P}$ in this section.
	\begin{definition}
		A commutative ring $R$ with unity is said to be a regular ring (Von-Neumann sense) if for every element $a\in R$, there exists an $x\in R$ such that $a = a^2x$.	
		
	\end{definition}
	A space $X$ is called $P$-space if every $G_\delta$-set in $X$ is open.
	\begin{definition}
		A $\tau \mathcal{P}$-space $(X,\tau, \mathcal{P})$ is called $\mathcal{P}P$-space if $C(X)_\mathcal{P}$ is regular.
	\end{definition}

	It has been observed in \cite{GJ1976} that for a Tychonoff space $X$, $X$ is a $P$-space if and only if $C(X)$ is regular. This fails when $X$ is not Tychonoff which can be seen in the following counter example.
	\begin{counter example} \label{6.01}
		Let $X=\mathbb{R}$ equipped with co-finite topology. Then $C(X)$ consists of all real valued constant functions on $\mathbb{R}$. Therefore, $C(X)$ is isomorphic to the field $\mathbb{R}$ which is regular. Thus $C(X)$ is regular. We define $G_r=X\setminus \{r\}$ for each $r\in \mathbb{Q}$. Then $G_r$ is open in $X$ for all $r\in \mathbb{Q}$, and $G=\bigcup_{r\in \mathbb{Q}}G_r=\mathbb{R} \setminus \mathbb{Q}$ is a $G_\delta$-set which is not open in $X$. Thus $X$ is not a $P$-space.
	\end{counter example}
	
	We note that Theorem 6.1 in \cite{GGT2018} also fails when $X$ is  not a Tychonoff space.
	To see this, we consider the next counter example.
	
	\begin{counter example} \label{6.02}
		Let $X=\mathbb{Q^*}=\mathbb{Q}\cup \{\infty \}$, the one-point compactification of $\mathbb{Q}$. Then every function in $C(\mathbb{Q^*})$ is constant. Thus $C(X)$ is isomorphic to $\mathbb{R}$, which is regular. However, $C(\mathbb{Q})$ is not regular, even though $\mathbb{Q}$ is a subspace of $\mathbb{Q^*}$. 
		Next, we show that $C(X)_F$ is not regular. \\We define $f(x)=\begin{cases}
			\sin (\pi x) &\text{ if }x\in \mathbb{Q} \\ 0 &\text{ if }x=\infty\end{cases}$. Then $f\in C(X)_F$. If possible, let there exists $g\in C(X)_F$ such that $f=f^2g$, then $g(x)=\frac{1}{\sin (\pi x)}$ for all $x\in X\setminus \mathbb{Z}$.
		For any $n\in \mathbb{Z}$, $g$ is unbounded in any neighbourhood of $n$ in $\mathbb{Q^*}$. Therefore, $D_g\supseteq \mathbb{Z}$, which contradicts $g\in C(X)_F$.
		This shows that the regularity of $C(X)$ might not imply the regularity of $C(X)_F$, that is, a $P$-space may not be an $\mathcal{F}P$-space.
	\end{counter example}

	However, if we assume $X$ to be Tychonoff, then the following is true.
	\begin{example}
		If $X$ is a $P$-space, then it is $\mathcal{P}P$-space.
		\begin{proof}
			Let $X$ be a $P$-space and $f\in C(X)_\mathcal{P}$. Then $f\in C(X\setminus \overline{D_f})$, where $X\setminus \overline{D_f}$ is a $P$-space, as it is a subspace of a $P$-space (by 4K in \cite{GJ1976}). So, $C(X\setminus \overline{D_f})$ is regular. Therefore, $\exists$ $g\in C(X\setminus \overline{D_f})$ such that $f|_{X\setminus \overline{D_f}}=(f|_{X\setminus \overline{D_f}})^2g$. Define $g^*$ on $X$ by $g^*(x)=\begin{cases}
			g(x) &\text{ when } x\in X\setminus \overline{D_f}
			\\	\frac{1}{f(x)} &\text{ when } x\in \overline{D_f}\setminus Z_\mathcal{P}(f)
			 \\ 0 &\text{ when } x\in \overline{D_f}\cap Z_\mathcal{P}(f)
			\end{cases}$. Then $D_{g^*}\subseteq \overline{D_f}$. Therefore, $g^*\in C(X)_\mathcal{P}$ and $f=f^2g^*$. Thus $C(X)_\mathcal{P}$ is regular, and so $X$ is a $\mathcal{P}P$-space.
		\end{proof}
	\end{example}
	The following result gives a generalisation of Theorem 6.2 (1)$\iff$(2) \cite{GGT2018}.
	\begin{theorem} \label{6.1}
		$X$ is a $\mathcal{P}P$-space if and only if for any zero set $Z$ in $X$, there exists a set $P\in \mathcal{P}$ such that $Z\setminus P$ is a clopen set in $X\setminus P$.
	\end{theorem}
	\begin{proof}
		Let $X$ be a $\mathcal{P}P$-space and $f\in C(X)_\mathcal{P}$. Then there exists $g \in C(X)_\mathcal{P}$ such that $f = f^2g$. Let $P = \overline{D_f}\cup \overline{D_g}\in \mathcal{P}$. Now $f|_{X\setminus P}, (\boldsymbol{1}-fg)|_{X\setminus P}$ are continuous. Also, $Z_\mathcal{P}(f)\setminus P=Z(f|_{X\setminus P})=(X\setminus P)\setminus Z((\boldsymbol{1}-fg)|_{X\setminus P})$ is a clopen subset of $X\setminus P$.
		Conversely, let the given condition hold and let $f\in C(X)_\mathcal{P}$. Then $Z_\mathcal{P}(f)\setminus P$ is a clopen subset of $X\setminus P$ for some $P\in \mathcal{P}$.
		Define $g \colon X \longrightarrow \mathbb{R}$ by
		$g(x)=\begin{cases}
			\frac{1}{f(x)} , &x\notin Z_\mathcal{P}(f) \\
			0, &x\in Z_\mathcal{P}(f)
		\end{cases}.$
		Then $\overline{D_g}\subseteq \overline{P}\cup \overline{D_f}\in \mathcal{P}\implies$ $\overline{D_g}\in \mathcal{P}$. So, $g\in C(X)_\mathcal{P}$ and $f = f^2g$. Hence $X$ is a $\mathcal{P}P$-space.
		
	\end{proof}

	\section{$\aleph_0$-Self injectiveness of $C(X)_\mathcal{P}$}
	
	In this section, we establish conditions under which $C(X)_\mathcal{P}$ is $\aleph_0$-self injective. In order to achieve this, we need the following definitions and theory.
	
\begin{definition} \cite{AHM2009}
		A ring $R$ is said to be $\aleph_0$-self injective if every module homomorphism $\phi \colon I\longrightarrow R$ can be extended to a module homomorphism $\widehat{\phi}\colon R\longrightarrow R$ where $I$ is a countably generated ideal of $R$.
\end{definition}
	A lattice-ordered vector space or vector lattice is a partially ordered vector space where the order structure forms a lattice.
	\begin{definition} \cite{HJ1962}
			
		An element $x$ of a vector lattice $X$ is called a weak order unit in $X$ if $x\geq 0$ and also for all $y\in X$, $\inf\{x,|y|\}=0$ implies $y=0$.

	\end{definition}
	\begin{definitions} \cite{HJ1962}
		By a lattice-ordered ring $(A,+,.,\vee,\wedge)$, we mean a lattice-ordered group that is a ring in which the product of positive elements is positive. If, in addition, $A$ is a (real) vector lattice, then $A$ is said to be a lattice-ordered algebra.
		
		A lattice-ordered ring $A$ is said to be Archimedean if, for each non-zero element $a\in A$, the set $\{na \colon n\in \mathbb{Z}\setminus\{0\} \}$ is unbounded.
		
		By a $\phi$-algebra $A$, we mean an Archimedean, lattice-ordered algebra over the real field $\mathbb{R}$ which has identity element $1$ that is a weak order unit in $A$.
		
		A $\phi$-algebra $A$ of real valued functions is said to be uniformly closed if it is closed under uniform convergence.
	\end{definitions}
	
	\begin{definitions} \cite{H1971}
		Let $A$ be a $\phi$-algebra. We denote $\mathcal{M}(A)$ as the compact space of maximal absolutely convex ring ideals of $A$ carrying Stone topology. Further, we denote $\mathcal{R}(A)$ to be the set of all real ideals of $A$.
	\end{definitions}
	
	\begin{definitions} \cite{AHM2009}
		For a subset $Q$ of a ring $R$, $Ann(Q)=\{r\in R\colon qr=0 \text{ for all }q\in Q \}$. A subset, $P$ of a ring $R$ is said to be orthogonal if the product of any two distinct elements of $P$ is zero.
		Suppose $P$ and $Q$ are disjoint subsets of $R$ whose union is an orthogonal subset of $R$. Then, an element $a\in R$ is said to separate $P$ from $Q$ if \begin{enumerate}
			\item $p^2a=p$ for all $p\in P$, and \item $a\in Ann (Q)$.
		\end{enumerate} 
	\end{definitions}
	
	We shall use Theorem 2.3 in \cite{H1971} to show that $C(X)_\mathcal{P}$ is isomorphic to an algebra of measurable functions.
	
	We assume $X$ to be a $\mathcal{P}P$-space and a $\tau \mathcal{PU}$-space, so that $C(X)_\mathcal{P}$ is regular and closed under uniform convergence. Further, for any $f\in C(X)_\mathcal{P}\setminus \{\boldsymbol{0}\}$, there exists $x\in X$ such that $f(x)\neq 0$. So, the set $\{\boldsymbol{n}f \colon n\in \mathbb{Z}\setminus \{\boldsymbol{0}\} \}$ is unbounded.
	We have already seen that $C(X)_\mathcal{P}$ is a lattice ordered group. It is also easy to see that it forms a real vector space and for any two positive elements $f,g\in C(X)_\mathcal{P}$, $fg$ is also positive. 
	Also, $C(X)_\mathcal{P}$ has the identity element $\boldsymbol{1}$ which is clearly a weak order unit.
	Thus, $C(X)_\mathcal{P}$ is a $\phi$-algebra which is closed under uniform convergence. That is, $C(X)_\mathcal{P}$ forms a uniformly closed $\phi$-algebra. 
	
	We have also seen that all maximal ideals of $C(X)_\mathcal{P}$ are $z_\mathcal{P}$-ideals which are in turn absolutely convex. Therefore, all maximal ideals of $C(X)_\mathcal{P}$ are in $Max(C(X)_\mathcal{P})$.
	
	Define for each $p\in X$, $M_p=\{f\in C(X)_\mathcal{P} \colon f(p)=0 \}$. Then, $C(X)_\mathcal{P}/M_p$ is isomorphic to $\mathbb{R}$, for each $p\in X$. Thus, $M_p$ is a real maximal ideal, for each $p\in X$ and is thus a member of $\mathcal{R}(C(X)_\mathcal{P})$. 
	
	\begin{theorem} \cite[Theorem 2.3]{H1971} \label{7.1}
		The following conditions on the $\phi$-algebra $A$ are equivalent.
		\\(a) $A$ is uniformly closed, regular, and $\bigcap \{M\colon M\in \mathcal{R}(A) \}=\{0\}$. \\
		(b) $A$ is isomorphic to an algebra of measurable functions.
	\end{theorem}
	
	We have \begin{equation*}
		\bigcap_{p\in X} M_p = \{\boldsymbol{0} \} \implies \bigcap \{M\colon M\in \mathcal{R}(C(X)_\mathcal{P}) \}=\{\boldsymbol{0}\}
	\end{equation*}
	From the above theorem, we have $C(X)_\mathcal{P}$ is isomorphic to an algebra of measurable functions. \\ Next we show that $\aleph_0$-self injectiveness of a ring is invariant under ring isomorphism.
	
	\begin{theorem} \label{7.2}
		If a ring $R$ is isomorphic to an $\aleph_0$-self injective reduced ring $S$, then $R$ is also $\aleph_0$-self injective.
	\end{theorem}
	\begin{proof}
		We shall use Theorem 2.2 of \cite{K1997} to prove the result. \\ Let $\psi \colon S \longrightarrow R$ be the given isomorphism. It is easy to see that, since $S$ is reduced, so is $\psi(S)=R$. \\ Further, by Theorem 2.2 of \cite{K1997}, $S$ is regular. Therefore, $R=\psi(S)$ is also a regular ring. \\ Let us now consider two disjoint subsets of $R$, $P$ and $Q$ whose union is a countable orthogonal subset of $R$. Then, $\psi^{-1}(P)$ and $\psi^{-1}(Q)$ are disjoint and their union is countable. Also, for any $s,s'\in \psi^{-1}(P)\cup \psi^{-1}(Q)$ with $s\neq s'$, $\psi(s), \psi(s')\in P\cup Q$ with $\psi(s)\neq \psi(s')$ (since $\psi$ is injective). As $P\cup Q$ is orthogonal, \begin{equation*}
			\psi(s)\psi(s')=0\implies \psi(ss')=0\implies ss'=0 \text{, since }\psi \text{ is injective.}
		\end{equation*} Thus, $\psi^{-1}(P)\cup \psi^{-1}(Q)$ is orthogonal.
		As $S$ is $\aleph_0$-self injective, by Theorem 2.2 of \cite{K1997}, there exists $a\in S$ that separates $\psi^{-1}(P)$ from $\psi^{-1}(Q)$. \\ We now show that $\psi(a)$ separates $P$ from $Q$. \begin{enumerate}
			\item Let $p\in P$, then $\psi^{-1}(p^2\psi(a))=(\psi^{-1}(p))^2a=\psi^{-1}(p)$, as $a\in S$ that separates $\psi^{-1}(P)$ from $\psi^{-1}(Q)$. It follows from the injectivity of $\psi^{-1}$ that $p^2\psi(a)=p$.
			Thus, $p^2\psi(a)=p$ for all $p\in P$.
			\item Let $q\in Q$, then $\psi^{-1}(q)\in \psi^{-1}(Q)$. As $a\in S$ that separates $\psi^{-1}(P)$ from $\psi^{-1}(Q)$, $\psi^{-1}(q)a=0$ which shows that $q\psi(a)=0$. 
			Thus, $\psi(a)\in Ann(Q)$.
		\end{enumerate} 
		Thus, we get an element in $R$ ($\psi(a)$) that separates $P$ from $Q$. 
		It follows from Theorem 2.2 of \cite{K1997} that $R$ is $\aleph_0$-self injective.
	\end{proof}
	Finally we use the above theory to establish the following theorem.
	\begin{theorem} \label{7.3}
		For a $\tau \mathcal{PU}$-space, the following conditions are equivalent. \\ (a) $X$ is a $\mathcal{P}P$-space. \\ (b) $C(X)_\mathcal{P}$ is isomorphic to an algebra of measurable functions. \\ (c) $C(X)_\mathcal{P}$ is $\aleph_0$-self injective.
	\end{theorem}
	\begin{proof}
		(a) $\implies$ (b) follows from the above discussions. (b) $\implies$ (c) can be seen from Theorem 7 of \cite{AHM2009} and the above theorem. Finally, (c) $\implies$ (a) follows directly from Theorem 2.2 of \cite{K1997}.
	\end{proof}
	
	We cannot omit the condition that $X$ is a $\tau \mathcal{PU}$-space. This can be seen from the following example.
	
	\begin{counter example}
		Let us consider \begin{equation*}
			X=\mathbb{N}\cup \bigcup_{k\in \mathbb{N}} \{\frac{1}{n}+k\colon n\in \mathbb{N} \} 
		\end{equation*} endowed with the subspace topology inherited from $\mathbb{R}_u$. Also let $\mathcal{P}=\mathcal{P}_f$, that is, the ideal of all finite subsets of $X$. Then, $C(X)_\mathcal{P}=C(X)_F$, which is not uniformly closed (by Theorem 2.9 in \cite{MBM2021}). Now, we consider the following subsets of $X$: \[A=\bigcup_{k\in \mathbb{N}}\{\frac{1}{2n}+k \colon n\in \mathbb{N} \} \text{ and }B=\bigcup_{k\in \mathbb{N}}\{\frac{1}{2n-1}+k \colon n\in \mathbb{N} \}.\]
		Define $P=\{\chi_{\{x \}}\colon x\in A\}$ and $Q=\{\chi_{\{x \}}\colon x\in B\}$ Then, $P$ and $Q$ are disjoint and $P\cup Q$ is countable and orthogonal. If there exists an $f\in \mathbb{R}^X$ that separates $P$ from $Q$, then $f(A)=\{1\}$ and $f(B)=\{0\}$. Thus, every point in $\mathbb{N}$ is a point of discontinuity of $f$. Therefore $f\notin C(X)_F$. This ensures from Theorem 2.2 of \cite{K1997} that $C(X)_F$ is not $\aleph_0$-self injective.
	\end{counter example}

	\bibliographystyle{plain}

\end{document}